\documentstyle[amstex,righttag]{amsart}
\newtheorem{thm}{Theorem}[section]
\newtheorem{lemma}{Lemma}[section]
\newtheorem{prop}{Proposition}[section]
\newtheorem{cor}{Corollary}[section]
\newtheorem{thmb}{Theorem} 
\newtheorem{defin}{Definition}[section]
\theoremstyle{definition}

\newtheorem{remark}{Remark} 

\def\squarebox#1{\hbox to #1{\hfill\vbox to #1{\vfill}}}

\def \Spc {\Sphere_c}
\def \phg {phg}

\def \cl {\operatorname{cl}}

\def \supp {\operatorname{supp}}

\def \sch {{\cal S}}

\def \mod {\text{mod}}

\def \S1 {\mbox{S}^{1}}

\def \om {\omega}

\def \omb {\overline{\omega}}
\def \omn {\omega_{n}}
\def \th {\theta}
\def \thn {\theta_{n}}
\def \thb {\overline{\theta}}

\def \Real {{\Bbb R}}
\def \Sphere {{\Bbb S} }
\def \Complex {{\Bbb C} }

\def \rest {|}
\def \Diffr {\operatorname{Diff}_r}

\def \pr {\prime}

\def \Natural{{\Bbb N}}

\renewcommand{\Re}{\operatorname{\rm Re}\nolimits}
\renewcommand{\Im}{\operatorname{\rm Im}\nolimits}

\def \y0l{Y_{0l}}
\def \dy0l{\Delta^{Y_{0l}}}
\def \ud {u_{\delta}}

\def \lzr{\langle z \rangle}
\def \lyr{\langle y \rangle}
\def \sca {\text{sc}}

\def \WF {\operatorname{WF}}

\def \dxt{\Delta_{\xi/|\xi|}}
\def \ep {\epsilon}

\title{Scattering on stratified media: the micro-local properties of the
scattering matrix and recovering asymptotics of perturbations}
\author{T. Christiansen and M. S. Joshi}
\begin{document}
  \maketitle

\begin{abstract}
The fixed energy scattering matrix is defined on a perturbed stratified
medium, and for a class of perturbations, its main part is shown to be
a Fourier integral operator on the sphere at infinity. This is facilitated by
developing a refined limiting absorption principle. The symbol of the
scattering matrix is shown to determine the asymptotics of a large class
of perturbations.
\end{abstract}

\section{Introduction}

In this paper, we study the structure of the scattering matrix on a perturbed
stratified medium. In particular, we show that its main part is a Fourier
integral
operator. En route to proving this theorem, we develop an improved limiting
absorption principle for a large class of perturbations, using techniques
of Fourier and microlocal analysis.
 As an application of our results, we prove that
the asymptotics of a perturbation can be recovered from the scattering matrix
at one
energy.

We recall that a stratified medium is a model space in which sounds waves
propagate
with a variable sound speed which depends on only one
coordinate. Thus, if we write the
coordinates on $\Real^{n}$ as $z=(x,y)$ with $x \in \Real^{n-1}$ and $y \in
\Real,$
we take the wave speed to be of the form $c_{0}(y)$ and study the wave equation
 \begin{equation}
(-\frac{\partial^2}{\partial t^2}-c_0^2\Delta )w=0,
\end{equation}
where $\Delta = -\sum_{j=1}^n
\frac{\partial^2}{\partial z_j^2.}$
We assume that $c_{0}$ is constant for $|y|$ large and that it is piecewise
smooth.
We let
\begin{align}
c_{+} = \lim \limits_{y \to \infty} c_{0}(y), \\
c_{-} = \lim \limits_{y \to -\infty} c_{0}(y).
\end{align}
In general, we do not require that $c_{+}$ be equal to $c_{-},$
but some of our results are stronger when they are equal.

A perturbed stratified medium is a medium on which the variable sound speed,
$c,$ has
the property that $c - c_{0}$ is well-behaved at infinity. Many previous papers
have
studied the case where the perturbation
$c-c_{0}$ is rapidly decaying. In particular,
precise asymptotics for $(c^2\Delta -(\lambda-i0)^2)^{-1}f$,
when $f\in {\cal S}(\Real^n)$, were proved in \cite{tsm}.
The inverse scattering problem for exponentially decaying perturbations was
studied
in \cite{bel,bbel,g-r,hi,we},
where it was shown that under certain conditions, the perturbation
can be recovered.

Here we study the case where the perturbation
$c-c_0$ has an asymptotic expansion in homogeneous terms
at infinity.
Under certain conditions on $c$ and
$c_0$ made more precise in Section \ref{s:aan}, we show that the
scattering matrix for $c^2\Delta$ is a Fourier integral operator and describe
its singular set.  Moreover, we show that the asymptotics of the perturbation
can be recovered from the scattering matrix at fixed energy.
We also establish the
lead term of the asymptotics for the limiting absorption principle.

 Our results use techniques developed by
Joshi and S\'a Barreto,
\cite{coulomb, smagailess, recpoten, metric, magnetic},
to study inverse problems in other settings, which build on work
by Melrose, \cite{sslaes}, and Melrose-Zworski, \cite{smagai}, on the structure
of the
scattering matrix on asymptotically Euclidean spaces.
 As in those inverse results, the
fundamental idea here is to compute
the symbol of the scattering matrix by solving transport equations
along geodesics on the sphere at infinity.  These equations express the
propagation of growth at infinity.

The analysis here is, however,
 considerably more involved as the unperturbed wave speed $c_0$
is not smooth on the compactified space achieved by adding the sphere at
infinity, even when
$c_0(y)$ is a smooth function of $y.$ This is because $c_0$ does not have nice
asymptotics
in $|z|.$ The upshot of this is that $c_0$ is well-behaved on the compactified
space
only after the space has been blown-up on the equator at infinity. This
manifests itself
in our analysis by requiring the geodesic flow at infinity to be refracted and
reflected by the equator. It was also seen in \cite{tsm} that it makes the
asymptotics in the limiting absorption principle much more complicated. There
is
a certain similarity here with many-body scattering,
compare, e.g. \cite{vasysr}. There the scattering problem is
complicated by the presence of a potential that does not decay in certain
directions
and thus appears as a spike on the sphere at infinity which causes refractions
and
reflections of the geodesic flows, \cite{vasysr}.
Indeed, the case where $c_{+}=c_{-}$ bears much resemblance to the many-body
case. However,
when $c_{+}$ and $c_{-}$ are different, there are effectively different energy
levels in the two hemi-spheres,
 which introduces new complications not present in the
many-body setting, and much of this paper is dedicated to coping with those
complications.

In Section \ref{ss:asm}
 we define the scattering matrix, and its ``main part.''  In
case the operator $c_0^2(D_y^2+\rho^2)$ has no eigenvalues as an operator
on $L^2(\Real, c_0^{-2}dy)$, then the main part of the scattering matrix
is the same as the scattering matrix.

Our first main result is
\begin{thm}\label{thm:sosm} Suppose $c,c_0$ satisfy the general
assumptions of Section \ref{s:aan}, and either hypothesis (H1) or (H2).
Then, if $c_+=c_-$, the main part
of the scattering matrix is a zeroth order Fourier integral operator associated
with
broken geodesic flow at time $\pi$.  If $c_->c_+$, then the main part
of the scattering matrix is a sum of Fourier integral operators associated
with the mapping $$(\omb, \omn)\mapsto (-\omb,\omn)$$ and the mapping
$$(\omb, \omn)\mapsto (-c_- \omb/c_+,-\sqrt{1-c_-^2|\omb|^2/c_+^2})$$ if
$\sqrt{1-c_+^2/c_-^2}<\omn$ and $$(\omb, \omn)\mapsto (-c_+ \omb /c_-,
\sqrt{1-c_+^2|\omb|^2/c_-^2})$$ if $\omn<0$.
\end{thm}
Here, when $c_+=c_-$, the geodesic flow is broken at the equator
$(\omb, 0)\subset \Sphere^{n-1}$.  This can be compared to
the situation for the Laplacian (\cite{smagai}),
or a perturbation of the
 Laplacian to an integral
 power (\cite{ch-j1}), on a manifold with asymptotically
Euclidean ends, where the scattering matrix is a zeroth order Fourier integral
operator
associated to geodesic flow at time $\pi$ on the the boundary ``at
infinity.''  An additional analogy is to $3$-body scattering, where the
three-cluster to three-cluster part of the scattering matrix is a sum of
Fourier integral operators associated to broken geodesic flow at time $\pi$
(\cite{vasysr}).  Other results on the structure of the scattering matrix
in $n$-body scattering may be found in \cite{vasymb}

Further results on the structure of the
 scattering matrix are given in
Proposition \ref{p:notmain}.

Our central inverse result is
\begin{thm}\label{thm:inverse}
Suppose $c$ and $c_0$ satisfy the general
assumptions of Section \ref{s:aan}, as well as either hypothesis (H1) or (H2),
and $n\geq 3$.  Then,
if $c_+=c_-$, the asymptotic expansion at infinity of $c-c_0$ is uniquely
 determined by $c_0$ and the transmitted singularities of the main part of the
scattering matrix at fixed nonzero
energy.  If $c_+<c_-$, then the asymptotic expansion is uniquely determined by
$c_0$ and the reflected singularities of the main part of the scattering matrix
at fixed nonzero energy.
\end{thm}
The reflected singularities are those associated
to the mapping $(\omb, \omn)\mapsto(-\omb, \omn)$ and, for $c_+=c_-$, the
transmitted singularities are those associated to the
mapping $\om \mapsto -\om$.
Corollary \ref{c:inverse} shows that knowledge of $c_+$, $c_-$, and
the singularities of the scattering matrix at fixed nonzero energy determine
$c$, within the class we consider, modulo a function vanishing faster than
the reciprocal of any polynomial at infinity.


Following the approach to studying the scattering matrix introduced in
\cite{smagai}, in Section \ref{s:cpo} we construct a parametrix
for the Poisson operator.  This is a key part of our proofs, as it facilitates
an understanding of the singularities of the scattering matrix.
We work particularly by adapting the techniques of \cite{smagailess} which are
essentially a concretization of the approach introduced in \cite{smagai}.
However, the different behaviour of the unperturbed operator
$ c_0^2\Delta$  in different regions at infinity means that
the analysis is considerably more involved.

To pass from a parametrix to the actual Poisson operator, we need a good
understanding of the behaviour
of $(\Delta -(\lambda-i0)^2c^{-2})^{-1}f$ at infinity, when $f\in
\lzr^{-\infty}L^2(\Real^n)$ and $(1-\phi(y))f\in {\cal S}(\Real^n)$
for some $\phi \in C_c^{\infty}(\Real)$.  In practice, the
$f$ for which we apply this will be
the error from the parametrix of the Poisson operator.  When
$c_+=c_-$, we can do this
by modifying some $n$-body results of \cite{g-i-s} and \cite{vasyab}.
However, when $c_+<c_-$ these results no longer apply, and we develop new
techniques.  The essential idea of these techniques is to repeatedly develop
better approximations with improving smoothness properties.
Thus Section \ref{s:bai} is devoted to understanding
$(\Delta -(\lambda-i0)^2c^{-2})^{-1}$, allowing us to finish the proof of
Theorem \ref{thm:sosm}. In particular, we prove the following limiting
absorption
principle
\begin{thm}\label{thm:smooth}\label{thm:smc}
  Let $c$ and $c_0$ satisfy the hypotheses of
Section \ref{s:aan} and hypothesis (H1) or (H2).
For any $\chi \in C_c^{\infty}(\Spc^{n-1})$, $f\in \lzr^{-\infty}L^2
(\Real^n)$,
$(1-\phi(y))f\in {\cal S}(\Real^n)$ for some $\phi \in C_c^{\infty}(\Real)$,
we have
$$\chi(z/|z|)(\Delta -c^{-2}(\lambda-i0)^2)^{-1}f
= e^{-i\lambda|z|/c}|z|^{-(n-1)/2}a_0(z/|z|) +u_1$$
where $a_0 \in C^{\infty}(\Spc^{n-1})$ and $u_1 \in \lzr^{\ep}L^2(\Real^n)$
for all $\ep >0$.
\end{thm}
Here $C_{c}^{\infty}( {\Bbb S}^{n-1}_c)$ is the space of smooth functions
vanishing
in a neighbourhood of the equator and in a neighbourhood of
$\{(\omb, \omn)\in \Sphere^{n-1}:\omn=\sqrt{1-c_+^2/c_-^2})\}$.

In Section \ref{s:inverse}, we use a modification
of some techniques of \cite{recpoten} to prove Theorem \ref{thm:inverse}
and apply some one-dimensional scattering theory to give further inverse
results.

An announcement of some of these results and an outline of part of
the proof can be found in the lecture notes \cite{ch-j2}.

We are grateful to Fritz Gesztesy for helpful conversations and
providing useful references and to Jim Ralston
for helpful discussions.
We thank the London Mathematical Society for supporting
this collaborative research through its small grants scheme.
The first author is grateful to the N.S.F. for partial support.

\section{Assumptions and Notation}\label{s:aan}
Throughout, $z=(x,y)\in \Real^{n-1} \times \Real$.

Both sound speeds $c$ and $c_0$ satisfy $0\leq c_m < c,c_0<c_M<\infty$.
Moreover, $c_0(y)$ is piecewise smooth and there exists a finite
$y_M$ so that $c_0(y)=c_{\pm}$ when $\pm y>y_M$, with $c_-\geq c_+$.
Moreover, all derivatives of $c_0$ are bounded except at finitely many
values of $y$.  This allows $c_0$ to be piecewise constant, for example.

We require that, away from the hypersurface
$\{y=0\}$, $c-c_0$ be smooth outside of a compact set $K$,
and for simplicity we choose $y_M$ so that
$K\subset \Real^{n-1}\times [-y_M,y_M]$.
Moreover, we make requirements on the behaviour of $c-c_0$
at infinity.  We have, for $y\not = 0$,
\begin{equation}\label{eq:catinfinity}
D_{z}^{\alpha}(c(z)-c_0(y))=
D_{z}^{\alpha}\sum_{j\geq J}^{N}\gamma_j(\frac{z}{|z|})|z|^{-j} +{\cal
O}(|z|^{-N-1-|\alpha|})
\end{equation}
for any $N$ and
any multiindex $\alpha$, where $\gamma_j
\in C^{\infty}_b(\Sphere^{n-1}\setminus
\{(\omb, 0)\})$.  Here we use the notation that
$C_b^{\infty}(X)$ is the space of smooth functions on $X$ that have all
derivatives bounded.
We shall take $J$ at least $2$ everywhere, although sometimes we shall
require it to be larger.  Some of our results hold under less restrictive
hypotheses.

Additionally, we shall often use one of the following hypotheses:\\
(H1) $J=2$, $c_+=c_-$, $c$ and $c_0$ are smooth.\\
(H2) $J\geq 4$.

We warn the reader that the choice of the total space dimension to
be $n$ rather than $n+1$ is in disagreement with \cite{tsm} and many other
papers on the subject.

We use the notation $\langle w \rangle =(1+|w|^2)^{1/2}$.  Throughout,
$\epsilon$ shall stand for a small positive quantity and $C$ for a positive
constant, either of which may change from line to line.

\section{Spectral theory of $c_0^2\Delta$}\label{s:st}

In order to define the (absolute) scattering matrix for $c^2\Delta$, we
will need some understanding of the generalized eigenfunctions of
$c_0^2\Delta $ and $c^2\Delta$, particularly of the space that parameterizes
them.  Further details can be found in, for example, \cite{b-d-g,d-g,
wederbook, wilcox}.

The operators $c_0^2\Delta$ and $c^2\Delta$ are formally self-adjoint on
$L^2(\Real^n,c_0^{-2}dz)$ and $L^2(\Real^n,c^{-2}dz)$, respectively and
have a unique self-adjoint extension.

Roughly speaking, the spectral measure of $c_0^2\Delta$ can be given
in terms of two kinds of families of functions.  At fixed
energy $\lambda$, the first is parameterized
by $\Sphere^{n-1}_c$, almost as for the Laplacian, though the
generalized eigenfunctions are more complicated.  Here
$\Spc^{n-1}=\{ \om=(\omb, \omn) \in \Sphere^{n-1}:\omn\not =0,\; \omn \not
= \sqrt{1-c_+^2/c_-^2}\}$.  (Compare \cite[Section 2.1]{wederbook}.)

A second type of generalized eigenfunction comes from eigenvalues of
$c_0^2(\kappa^2+D_y^2)$ on $L^2(\Real, c^{-2}_0dy)$,
if there are any.  If there are any eigenvalues, let
$\lambda_1^2(\kappa)<\lambda_2^2(\kappa)
<\cdot\cdot\cdot < \lambda_{k(\kappa)}^2(\kappa)<c_+^2\kappa^2$
denote the eigenvalues of $c_0^2(\kappa^2+D_y^2)$.  There may not
be any eigenvalues, but if there are, there are only finitely many for
fixed $\kappa$ and the number grows with $\kappa^2$.  Additionally,
if $\kappa>0$ and $\lambda_j>0$, then
$\frac{d \lambda_j}{d\kappa}>0$, as can be seen by an
integration by parts argument (see, e.g., \cite[Sect. 2.2]{tsm}).

Let $\kappa_j^0$ be the smallest positive number such that
$c_0^2(\kappa^2+D^2_y)$ has $j$ eigenvalues for all $\kappa >\kappa_j^0$.
Let $\kappa_j$ be the inverse of $\lambda_j$ (with the same sign), and let
$t_j=
\lim_{\kappa\downarrow \kappa_j^0}\lambda_j^2(\kappa)=c_+^2(\kappa_j^0)^2$.
The $\{t_j\}$ are called thresholds of $c_0^2\Delta$.  Let
$T(\lambda)$ be the number of thresholds $t_j$ less than $\lambda^2$.
For $0<j\leq T(\lambda)$, $\omb\in \Sphere^{n-1}$, let
$$\Phi_j(z,\lambda,\om)=e^{i\kappa_j(\lambda)x\cdot \omb}f_j(y),$$
where $f_j(y)\in L^2(\Real, c_0^{-2}dy)$ satisfies
$$c_0^2(\kappa_j^2+D_y^2)f_j=\lambda^2 f_j$$
and note that $(c_0^2\Delta -\lambda^2)\Phi_j=0$.

At energy level $\lambda^2$, we can parameterize the generalized
eigenfunctions of $c_0^2\Delta$ by $\Spc^{n-1}$ and $T(\lambda)$
copies of $\Sphere^{n-2}$.
The continuous spectrum of $c^2\Delta$ is parameterized by the same
space as that of $c^2_0\Delta$.

Because of the described parametrization of the continuous spectrum at
fixed energy, the (absolute) scattering matrices of $c_0^2\Delta $
and $c^2\Delta$ are operators from
$L^2(S_c)\oplus_{1\leq j \leq T(\lambda)}L^2(\Sphere^{n-2})$ into itself.
In \cite{tsm} a definition of the scattering matrix is given in terms of
the generalized eigenfunctions.  Here, however, it will be more useful to
define the (absolute) scattering matrix using the Poisson operator,
which we shall do in Section \ref{ss:asm}.

\section{The Poisson operator and the scattering matrix}

The Poisson operator is defined as an operator
$$P(\lambda):C_c^{\infty}(\Spc^{n-1})
\oplus _{i=1}^{T(\lambda)}C^{\infty}(\Sphere^{n-2})\rightarrow \langle z
\rangle ^{-1/2-\epsilon}L^2(\Real^n).$$
\begin{defin}\label{d:po}
If $g=(g_0,g_1,...,g_{T(\lambda)})\in
C_c^{\infty}(\Spc)
\oplus _{i=1}^{T(\lambda)}C^{\infty}(\Sphere^{n-2}),$ then $P(\lambda)g=u$,
$(c^2\Delta -\lambda^2)u=0$,
and $u$ is determined by its asymptotics at infinity:
$$u \sim |z|^{-(n-1)/2}e^{i\lambda |z|/c}g_0\left(\frac{z}{|z|}\right)
+|x|^{-(n-2)/2}\sum_{j=1}^{T(\lambda)}e^{i\kappa_j(\lambda)|x|}f_j(y)g_j
\left( \frac{x}{|x|}\right)+u_0+\sum_{j=1}^{T(\lambda)}u_j +\tilde{u}.$$
The functions $u_0,$ $u_j$, $\tilde{u}$ satisfy
\begin{equation*}\begin{array}{c}
(\frac{\partial}{\partial |z|}+i\lambda/c)u_0 \in \lzr^{\ep}L^2(\Real^n);\;
\; u_0 \in \lyr^{1/2+\ep}\lzr^{\ep}L^2(\Real^n);\\
u_j = |x|^{-(n-2)/2}e^{-i\kappa_j(\lambda)|x|}b_j(x/|x|)f_j(y),\; 1 \leq j \leq
T(\lambda);\\
\tilde{u},\frac{\partial}{\partial |z|} \tilde{u} \in \lzr^{\ep}L^2(\Real^n),
\end{array}\end{equation*}
for any $\ep >0$.  Here $\kappa_j$, $f_j$ are as in Section \ref{s:st}.
\end{defin}
Proposition \ref{p:powd} shows that this expansion
 uniquely determines the Poisson
operator.  Definition \ref{d:sm}, using Proposition \ref{p:npoe},
defines the (absolute) scattering matrix via the Poisson operator.
The existence of the Poisson operator is proved in Section \ref{s:bai}.

\subsection{The Poisson operator is uniquely determined}

In order to show that the Poisson operator above is indeed well-defined,
we shall need a uniqueness result, for whose proof we shall use Proposition
\ref{p:ep}.

In proving the following proposition, we shall use some results
of Weder, \cite{wederlapat,wederbook} (See also \cite{db-pii}.).
We recall some of his results below.

Let $A=(-i/4)(z\cdot  \nabla_z + \nabla _z \cdot z)$.  We define the
commutator $[\Delta-\lambda^2/c^2,A]$ as a quadratic form (See the proof of
Theorem 5.4, \cite{wederbook}.).
By \cite[Lemma 3.1]{wederlapat} for all $\lambda >0$, $\mu>-\lambda^2/c_-^2$,
 there is a compact operator $K$, a compact interval $\Lambda$
containing $\mu$, and $\beta >0$ such that
\begin{equation}\label{eq:me}
iE_{\Lambda}[\Delta -\lambda^2c^{-2},A]E_{\Lambda} \geq \beta E_{\Lambda}+K
\end{equation}
where $E_{\Lambda}=E_{\Lambda}(\Delta -\lambda^2c^{-2})$ is the
spectral projector for $\Delta -\lambda^2c^{-2}$.

The following proposition and its proof, included for the
convenience of the reader, are essentially
adapted from \cite[Lemma 4.17]{bbel}.
\begin{prop}\label{p:ep}
If $u\in \lzr^{\ep}L^2(\Real^n)$ for
 every $\ep >0$ and $(\Delta -\lambda^2/c^{-2})u=0$, then $u\equiv 0$.
\end{prop}
\begin{pf}
By the results of \cite{wederlapat, wederbook},
there is no nontrivial $L^2$ null space of
$\Delta -\lambda^2/c^2$, so it suffices to show that $u \in L^2(\Real^n)$.

For $\ep, \delta>0$, let $u_{\delta}= (1+\delta\lzr)^{-\ep}u\in L^2(\Real^n)$.
Let $L=\Delta -\lambda^2/c^2$, and let $\Phi \in C_c^{\infty}(\Real)$
be $1$ in a neighbourhood of $0$.

Note that
$$L\Phi(L)\ud = \Phi(L)\left( \sum _j \frac{2\ep \delta z_j}{\lzr}
(1+\delta \lzr)^{-1}\frac{\partial }{\partial z_j}\ud + f_{\ep \delta }(z)u
\right)$$
where
\begin{equation}\label{eq:fest}
|f_{\ep \delta }(z)|\leq C\frac{\ep}{\lzr^2}(1+\delta \lzr)^{-\ep}
\end{equation}
and the constant $C$ is independent of $\ep$ and $\delta $.

Then
\begin{multline}
([L,A]\Phi(L)\ud, \Phi(L)\ud)=-
2i\Im (A L \Phi(L)\ud, \Phi(L)\ud)\\
= -2i\Im (A \Phi(L)( \sum _j \frac{2\ep \delta z_j}{\lzr}
(1+\delta \lzr)^{-1}\frac{\partial }{\partial z_j}\ud + f_{\ep \delta} (z)u
),\Phi(L)\ud).\end{multline}
Using this equality, (\ref{eq:fest}), and the fact that
$$\Phi(L):\lzr^{-\gamma}L^2(\Real^n)\rightarrow \lzr ^{-\gamma}L^2(\Real^n),$$
we obtain
\begin{equation}
\label{eq:ucb}
|([L,A]\Phi(L)\ud, \Phi(L)\ud)| \leq \ep C_{\Phi}\|\ud\|\|\Phi(L)\ud\| +C_2.
\end{equation}
Here and below
 $C_{\Phi}$, $C_2$ are constants which
may change, independent of $\ep$ and $\delta$,
but depending on $\Phi$, and $C_2$ depends on $\| \lzr^{-\ep_0}u\|$ as well.
Since $(1-\Phi(L))u=0$, we have
\begin{equation}\label{eq:uude}
\|\ud\| = \| \Phi(L)\ud + [\Phi(L),(1+\delta \lzr)^{-\ep}] u \|
\leq \| \Phi(L)\ud\| +C_2.
\end{equation}

However, using (\ref{eq:me}) and the fact that $0$ is not an eigenvalue
of $L$, we obtain
$$([L,A]\Phi(L)\ud,\Phi(L)\ud)\geq \beta_1(\Phi(L)\ud,\Phi(L)\ud)$$
for some $\beta_1 >0$, if the support of $\Phi$ is chosen sufficiently small.

By choosing $\ep$ sufficiently small, then, and using (\ref{eq:ucb})
and (\ref{eq:uude}), we get that
$$\beta_1 \| \Phi(L)\ud\| ^2 \leq \frac{\beta_1}{2} \| \ud\|^2 + C_{2}.$$
Then $$\| \Phi(L)\ud\| \leq C_{2},$$
and using (\ref{eq:uude}), this shows that for
sufficiently small $\ep$, $\| \ud \|$ is bounded by a constant
independent of $\delta$, and thus $u\in L^2(\Real^n)$, and $u \equiv 0$.
\end{pf}

We shall use the following notion of an ``outgoing'' function.
\begin{defin}\label{d:outgoing}
 A function $u\in \lzr^{1/2+\ep}L^2(\Real^n)$ will be called
outgoing if it has a decomposition
$u=u_0+\sum_{1}^{T(\lambda)}u_j+\tilde{u}$ with the following properties
\begin{equation*}\begin{array}{c}
(\frac{\partial}{\partial |z|}+i\lambda/c)u_0 \in \lzr^{\ep}L^2(\Real^n);\;
\; u_0 \in \lyr^{1/2+\ep}\lzr^{\ep}L^2(\Real^n);\\
u_j = |x|^{-(n-2)/2}e^{-i\kappa_j(\lambda)|x|}b_j(x/|x|)f_j(y),\; 1 \leq j \leq
T(\lambda);\\
\tilde{u},\frac{\partial}{\partial |z|} \tilde{u} \in \lzr^{\ep}L^2(\Real^n),
\end{array}\end{equation*}
for any $\ep >0$.
\end{defin}

\begin{prop}\label{p:powd}  Given $f\in L^2(\Real^n)$,
there is at most one outgoing $u\in \lzr^{1/2+\ep}L^2(\Real^n)$ with
$(\Delta -\lambda^2/c_0^2)u=f.$
\end{prop}
\begin{pf}
Suppose there are two such $u$.  Then by considering the difference we
can reduce this to the case $f \equiv 0$.  Then
\begin{multline}
0 =\int _{|z|<R}(\Delta -\lambda^2/c^2)u \overline{u}= -
\int_{|z|=R} (\frac{\partial}{\partial |z|}u \overline{u}-u
\frac{\partial}{\partial |z|} \overline{u })\\=
2\int _{|z|=R} \frac{i\lambda}{c}|u_0|^2+\sum_{j,k}i\kappa_j(\lambda)
u_j\overline{u}_k + i\Re\sum_j (\lambda/c+\kappa_j)u_0\overline{u}_j
+i\Im(u_0e_0 + \sum_{j}u_je_j) +i\Im \tilde{e}\tilde{f}
\end{multline}
where $e_0, e_j, \tilde{e}, \tilde{f} \in \lzr^{\ep}L^2(\Real^n)$ for
all $\ep >0$.
Using the facts that $\int f_j (y)\overline{f}_k(y)dy =0$ if $j\not = k$
and $f_j$ is exponentially decreasing,
this implies that as $R\rightarrow \infty$
\begin{equation}\label{eq:uineq}
\int _{|z|=R} \frac{2i\lambda}{c}|u_0|^2+\sum_{j}2i\kappa_j(\lambda)
|u_j|^2 \leq C \int_{|z|=R}(\sum_j |u_0 \overline{u}_j|
+|u_0||e_0| + \sum_{j}|u_je_j| +|\tilde{e}||\tilde{f}|) + {\cal O}(R^{-2-\ep}).
\end{equation}
Since $u_j \in \lyr^{-\infty}\langle x \rangle ^{1/2+\ep}L^2(\Real^n)$
and $u_0 \in \lyr^{1/2+\ep}\lzr^{\ep}L^2(\Real^n)$, we have
$u_0\overline{u}_j \in \lyr^{-\infty}\lzr^{1/2+2\ep}L^1(\Real^n)$.  Therefore,
the right hand side of (\ref{eq:uineq}), considered as a function
of $R$ for large $R$, is in $R^{1/2+2\ep}L^1(\Real_+)$, so that
$u_0,\; u_j \in \lzr^{1/4+\ep}L^2(\Real^n)$.  This in turn means that
$b_j=0=u_j$ for $1\leq j \leq T(\lambda)$.  Now we have
$\frac{\partial}{\partial |z|}u = (-i\lambda /c) u +\tilde{u}$,
$\tilde{u } \in \lzr^{\ep}L^2(\Real^n)$, any $\ep >0$.
Repeating the argument almost as above, we obtain
$$\int_{|z|=R}|u|^2 \leq C\int_{|z|=R}|u\tilde{u}|
\leq C (\int_{|z|=R}|u|^2)^{1/2}(\int_{|z|=R}|\tilde{u}|^2)^{1/2}$$
and thus
$$(\int_{|z|=R}|u|^2)^{1/2}\leq C (\int_{|z|=R}|\tilde{u}|^2)^{1/2}
$$ so that $u\in \lzr^{\ep}L^2(\Real^n)$, any $\ep >0$.  By the previous
Proposition, $u\equiv 0$.
\end{pf}

\subsection{The Absolute Scattering Matrix}\label{ss:asm}
In order to use the Poisson operator to define the scattering matrix, we
shall need the following proposition, whose proof is a corollary of our
construction of Section \ref{s:cpo} and Theorem \ref{thm:smc} (see
Section \ref{s:smooththm}).
\begin{prop} \label{p:npoe}
 Let $g=(g_0,g_1,...,g_{T(\lambda)}) \in
C_c^{\infty}(\Spc^{n-1})\oplus_1^{T(\lambda)}C^{\infty}
(\Sphere^{n-2})$.
Let $\theta =\frac{z}{|z|}\in K$, for $K$ a compact set in $\Spc^{n-1}$.
Then, for $\theta \in K$, as $|z|\rightarrow \infty$,
$$P(\lambda)g = |z|^{-(n-1)/2}[e^{i\lambda |z|/c}g_0(\theta)_{\rest K}
+e^{-i\lambda |z|/c}(g_0'(\theta))_{\rest K}] + \tilde{u}_K$$
where $\tilde{u}_K\in \lzr^{\ep}L^2(\Real^n\cap(K\times [1,\infty)))$
for any $\ep >0$.\\
Let $y\in K_y$, $K_y\subset \Real$ compact, and let $\thb=\frac{x}{|x|}$.
Then, as $|x|\rightarrow \infty$,
$$P(\lambda)g=|x|^{-(n-2)/2}[\sum_1^{T(\lambda)}e^{i\kappa_j(\lambda)|x|}
g_j(\theta)f_j(y)_{\rest K_y}+ \sum_1^{T(\lambda)}e^{-i\kappa_j(\lambda)|x|}
g_j'(\theta)f_j(y)_{\rest K_y}]+ \tilde{u}_c$$
where $\tilde{u}_c\in \langle x \rangle^{\ep} L^2((\Real^{n}\cap(\{|x|>1\})
\times K_y))$ for any $\ep >0$.
\end{prop}

This information about the Poisson operator allows us to define the (absolute)
scattering matrix $A(\lambda)$.
\begin{defin}\label{d:sm}
 The (absolute) scattering matrix $A(\lambda)$ is given, for
$g\in C_c^{\infty}(\Spc^{n-1})\oplus_1^{T(\lambda)}C^{\infty}(S^{n-2})$,
by
$A(\lambda)g=g'\in L^{2}(\Spc^{n-1})\oplus L^2(S^{n-2})$,
where for any compact set $K\subset \Spc^{n-1}$, $(g'_0)_{\rest K}$
is as in Proposition \ref{p:npoe}, and, for $1\leq j \leq T(\lambda)$,
 $g'_j$ is as in Proposition
\ref{p:npoe}.
\end{defin}
We remark that this definition differs slightly from the absolute scattering
matrix discussed in \cite{tsm}.  However, as
the two differ by a straight-forward
normalization, we shall use this definition here both to emphasize the
similarities with the absolute scattering matrix as defined in \cite{sslaes}
and because it is most convenient for the inverse results.

For fixed $\lambda$, $A(\lambda)=(A_{ij}(\lambda))$,
$0\leq i,j\leq T(\lambda)$,
with the $A_{ij}$ operators.
We call $A_{00}(\lambda)$ the ``main part'' of
the scattering matrix.  If the operator $c_0^2(D_y+\rho^2)$ has no
eigenvalues on $L^2(\Real, c^{-2}_0dy)$, the ``main part'' of the scattering
matrix is just the scattering matrix.

\section{The approximate Poisson operator}\label{s:cpo}

For $g=(g_0,g_1,...,g_{T(\lambda)})\in L^2(\Spc^{n-1})\oplus_1^{T(\lambda)}
L^2(\Sphere^{n-2})$, let $\Pi_jg=g_j$.  Since our inverse results
involve the main part of the scattering matrix, $A_{00}(\lambda)$, we are
most interested in $P_0(\lambda)=P(\lambda)\Pi_0$, where
$P(\lambda)$ is the Poisson operator.  Here we construct an approximation
$\tilde{P}_0$ of $P_0$.

Let $z\in \Real^n$, $\lambda \in \Real$, $\lambda \not = 0$, and let
$\omega \in \Sphere^{n-1}$, $\om=(\omb, \omn)$,
$\omn\not = 0$.  We will construct a (partial) approximate
Poisson operator $\tilde{P}_0(z,\lambda, \om)$ such that
$$(c^2\Delta -\lambda^2)\tilde{P}\in \langle z \rangle ^{-\infty}L^2
(\Real^n),\;
(1-\chi(y))(c^2\Delta -\lambda^2)\tilde{P}_0\in {\cal S}(\Real^n)$$ where
$\chi(y)\in C_c^{\infty}(\Real)$ is $1$ for $|y|\leq y_M+1$,
and, so that distributionally as $|z|\rightarrow \infty,$
$$\tilde{P}_0(z,\lambda,\om)=|z|^{-(n-1)/2}
\left(e^{i\lambda |z|/c_0(y)}\delta_{\om}(\frac{z}{|z|})
+e^{-i\lambda|z|/c_0(y)}g(\om, \frac{z}{|z|})\right)+{\cal
O}(|z|^{-(n+1)/2}).$$

We will show how to construct such an approximation when $\omn>0$; the
case of $\omn<0$ is quite similar.  The construction involves solving
away errors at infinity.  Since the model operator $c_0^2\Delta$ has different
behaviour depending on the region ``at infinity'' ($y>y_M$, $|y|<y_M$, or
$y<-y_M$) the techniques involved necessarily depend on the region
in which $z$ lies.

We shall need the following
\begin{lemma}\label{l:eu}
For $f\in L^2([-y_M,y_M])$, the boundary value problem
$$c_0^2(\lambda^2|\omb|^2/c_+^2+D^2_y)b-\lambda^2b=f$$
with boundary conditions
\begin{align*}
-i\lambda \omn/c_+b(y_M)-b^{\pr}(y_M)&=\alpha_1\\
i\lambda (1/c_-^2-|\omb|^21/c_+^2)^{1/2}b(-y_M)-b^{\pr}(-y_M)&=\alpha_2
\end{align*}
has a unique solution in $L^{2}([0,y_0],c_0^{-2}dy)$ if
$1/c_-^2\geq |\omb|^2/c_+^2$.
\end{lemma}
\begin{pf}
This boundary problem can be reduced to the form
\begin{align*}
-i\lambda \omn/c_+b(y_M)-b^{\pr}(y_M)&=0\\
i\lambda (1/c_-^2-|\omb|^2/c_+^2)^{1/2}b(-y_M)-b^{\pr}(-y_M)&=0\\
c_0^2(\lambda^2|\omb|^2/c_+^2+D^2_j)b-\lambda^2b&=g.
\end{align*}
This has a solution if the adjoint operator has no nontrivial null space; and
the solution is unique if the homogeneous equation has no nontrivial
solutions.

The adjoint operator is
the operator
$$c_0^2(\lambda^2|\omb|^2/c_+^2+D^2_j)-\lambda^2$$
with domain
\begin{multline}
\{ h\in L^2([0,y_0],c_0^{-2}dy):
-i\lambda \omn/c_+h(y_M)+h^{\pr}(y_M)=0\\
\mbox{ and } i\lambda (1/c_-^2-|\omb|^2/c_+^2)^{1/2}h(-y_M)
+h^{\pr}(-y_M)=0
\}.
\end{multline}
Suppose $g$ is a nontrivial element of the null space of the adjoint operator.
Then
\begin{align}\label{eq:ibp}
0&=\int_{-y_M}^{y_M}(c_0^2(\lambda^2|\omb|^2/c_+^2+D^2_y)-\lambda^2)g
\overline{g}c_0^{-2}dy \nonumber \\
&=-g^{\pr}(y_M)\overline{g}(y_M)+g^{\pr}(-y_M)\overline{g}(-y_M)
+g(y_M)\overline{g}^{\pr}(y_M)-g(-y_M)\overline{g}^{\pr}(-y_M)
\nonumber
\\ &
+ \int_{-y_M}^{y_M}\left[(\lambda^2|\omb|^2/c_+^2-\lambda^2 c_0^{-2})|g|^2
+|D_yg|^2 \right]dy
\end{align}
Using the boundary conditions, we find that the second line is equal to
$$-2i\lambda \omn/c_+|g(y_M)|^2-2i\lambda(1/c_-^2-|\omb|^2/c_+^2)^{1/2}
|g(-y_M)|^2.$$
Since the third line of (\ref{eq:ibp})
is real, this means that $0=g(y_M)=g^{\pr}(y_M)$, or
$g(y)\equiv 0$.

A similar calculation shows that the original operator has no nontrivial
null space.
\end{pf}

Let $\om=(\omb,\omn)$ with
$\omn>0$.  In our construction of the approximation to the Poisson
operator $P_0$,
we begin with the function $\Phi(z,\lambda,\om)$ which is
defined by
$$
e^{i\lambda x\cdot \omb/c_+}\phi_+(y)$$ where
$\phi_+$ satisfies
\begin{equation}\label{eq:ode}
c_0^2(D^2_y+\lambda^2(1-\omn^2)c_+^{-2})\phi_+=\lambda^2\phi_+
\end{equation}
and, as $y\rightarrow \infty$,
\begin{equation}
\phi_+(y)\sim e^{i\lambda y\omn/c_+}+R_+(\lambda, \omn)e^{-i\lambda\omn y/c_+}
\end{equation}
and as $y\rightarrow -\infty$,
\begin{equation}\label{eq:yni}
\phi_+(y)\sim T_+(\lambda,\omn)
e^{i\lambda y\sqrt{1/c_-^2-1/c_+^2+\omn^2/c_+^2}}
\end{equation}
where when $1/c_-^2-1/c_+^2+\omn^2/c_+^2<0$ we take the square
root so that the right hand side of
(\ref{eq:yni}) is exponentially decreasing.  We have
$(c_0^2\Delta-\lambda^2)\Phi=0$.  Note that, up to a constant multiple
which depends only on $n$,
$\lambda$, and $c_{\pm}$, $\Phi$ is the Schwartz kernel of the (partial)
Poisson operator $P_{0,0}$
for $c_0^2\Delta$ when $\om$ is in the upper hemisphere of $\Sphere^{n-1}$.
We use this as our starting point.

When $y>y_M$,
we use the
techniques of \cite{smagailess} to construct $\tilde{P}$.  Note that when
we apply $c^2\Delta-\lambda^2$ to $\Phi$ we obtain an error which,
for $y>y_M$, is of the form
\begin{equation}\label{eq:e1t}
a_1 e^{i\lambda z\cdot \om/c_+} + a_2 e^{i\lambda (x\cdot \omb -y\omn)/c_+}
\end{equation}
where $a_1,a_2$ are classical polyhomogeneous symbols (in $|z|$)
of order $-2$.
We use the techniques of \cite{smagailess} to find a term of the form
$ae^{i\lambda z\cdot \om/c_+} $ with
$a\in S^{-1}_{phg}$ which will solve away the first term
in the error (\ref{eq:e1t}) here.  Just as in \cite{smagailess}, this
is done iteratively, solving away an error  $d\in e^{i\lambda z\cdot \om/c_+}
(S^{-j}_{phg}\mod S^{-j-1}_{phg})$ by using
$b e^{i\lambda z\cdot \om/c_+}$ with $b\in S^{-j+1}_{phg}$,
and $b$ solving
the transport equation
\begin{equation}\label{eq:te1}
-2ic_+\lambda\om\cdot \frac{\partial b}{\partial z}=d
\end{equation}
along the geodesics on the unit sphere at infinity.  We choose
$b$ so that $b$ is smooth at $z/|z|=\om$ in order to
keep the right coefficient of $e^{i\lambda |z|/c_+}$
in the distributional asymptotic expansion.

Let $s$ be the geodesic distance on the sphere from $\om$ and
let $\tilde{\theta}$ be angular coordinates about $\om$.
The equation (\ref{eq:te1}) can be solved, modulo
an error of lower order, by $b_{I,-j+1}|z|^{-j+1}$
 just as in \cite[Section 2]{smagailess},
giving
\begin{equation}\label{eq:ste}
b_{I,-j+1}(s,\tilde{\theta},\om)=\frac{i}{2\lambda c_+ (\sin s)^{j-1}}
\int_0^s (\sin s')^{j-2}d_{-j,I}(s',\tilde{\th};\om)ds'
\end{equation}
where $d_{-j,I}|z|^{-j}$ is the term of homogeneity $-j$ in the error.
Note that as long as
$z/|z|$ is in the upper
hemisphere we are away from $-\omega$ so the transport equations have
a smooth solution.

We find $b_{I,-j}$ iteratively and then use Borel's lemma to asymptotically
sum them, obtaining
a $b_{I}$ such that
$$(c^2\Delta -\lambda^2)(\Phi + b_Ie^{i\lambda z \cdot \om/c_+})
=e^{i\lambda x\cdot \omb/c_+}e^{-i\lambda y\omn/c_+} a_2 +
{\cal O}(|z|^{-\infty})$$
when $y>y_M$. Note that the construction of $b_I$ has not changed $a_2$.

We will apply almost the same technique to solve away the error
$e^{i\lambda x\cdot \omb/c_+}e^{-i\lambda y\omn/c_+}a_2$, $y>y_M$, away
from $z/|z|=(-\omb, \omn)$.  Here we will use solutions to the transport
equation where we choose the initial condition at $y/|z|=0$, and the
solutions, in analogy to (\ref{eq:ste}), are of the form
\begin{equation}
b_{R,-j+1}(s,\tilde{\th},\om)=
\frac{i}{2\lambda c_+ (\sin s)^{j-1}}
\left[ \int_{s_{R_0}}^s
(\sin s')^{j-2}d_{-j,R}(s',\tilde{\theta};\om)ds' +C_{R,j-1}
\right]. \end{equation}
Here $s$ is the distance on $\Sphere^{n-1}$ from $\theta=z/|z|$ to the point
$(\omb, -\omn)$ and $\tilde{\theta}$ is the angular coordinate about
$(\omb, -\omn)$.  The value $s=s_{R_0}$ corresponds to $\thn=0$, and
$C_{R,j}$ depends only on $\om$ and $\tilde{\theta}$. We postpone to
Section \ref{ss:poas}
discussion
of the form of the parametrix near $z/|z|=(-\omb, \omn)$.

In the lower hemisphere, we use a similar technique
if $1/c_-^2-|\omb|^2/c_+^2>0$.  Here the error
term is of the form
$$e^{i\lambda x \cdot \omb/c_+}e^{i\lambda
\sqrt{1/c_-^2-|\omb|^2/c_+^2}y}a_T,$$
where $a_T\in S^{-2}_{\phg}.$  Again we have solutions like (\ref{eq:ste})
to the transport equation, although this time $s$ measures the distance
on the sphere from the point $(c_-\omb/c_+, \sqrt{1-c_-^2|\omb|^2/c_+^2}).$
We will have a solution to
the transport equation away from $(-c_-\omb/c_+,
-\sqrt{1-c_-^2|\omb|^2/c_+^2})$
  of the form
$$e^{i\lambda x \cdot \omb/c_+}e^{i\lambda \sqrt{1/c_-^2-|\omb|^2/c_+^2}y}
\sum b_{T,-j}|z|^{-j},$$
where
$$
b_{T,-j+1}(s,\tilde{\theta},\om)=\frac{i}{2\lambda c_- (\sin s)^{j-1}}
\left[ \int_{s_{T_0}}^s (\sin s')^{j-2}d_{-j,T}(s',\tilde{\th};\om)ds'
+C_{T,j-1}\right].$$

The constants (in $s$) $C_{R,j}$ and $C_{T,j}$ are to be determined.  Of
course their values affect subsequent errors and thus subsequent $b_{R,-j}$,
$b_{T,-j}$.

We will use a different technique to construct the
solutions when $|y|<y_M$.  We choose the solutions so that $\tilde{P}_0$ is
$C^1$ on
$\Real^n$.
We point out that if $c_0$ is not smooth, for example, if it
is piecewise constant, we should expect it to be impossible to find a smooth
Poisson operator on $\Real^n$.

The values of $C_{R,j}$, $C_{T,j}$ are determined by solutions to
boundary value problems that arise in constructing the parametrix when
$|y|<y_M$, as described below.

When $|y|<y_M$, we look for an approximate solution of the
form \begin{equation}\label{eq:mid}
e^{i\lambda x\cdot \omb/c_+} \sum_{j\geq 0}
|z|^{-j}\left(
b_{M,-j}(\frac{x}{|x|},y)+|z|^{-1}\tilde{b}_{M,-j}(\frac{x}{|x|},y)\right).
\end{equation}
The term $|z|^{-1}\tilde{b}_{M,-j}$ is of lower order and is included
to improve the regularity at $y=\pm y_M$.
Note that
$$(c^2\Delta -\lambda^2)(|z|^{-j}b(\frac{x}{|x|},y)e^{i\lambda x\cdot
\omb/c_+})
= \left(c_0^2(D^2_y+\lambda^2|\omb|^2/c_+^2)b-\lambda^2 b\right)
e^{i\lambda x\cdot \omb/c_+}|z|^{-j}+{\cal O}(|z|^{-j-1}).$$
Therefore, for $|y|<y_M$, to solve away an error
of the form
$|z|^{-j}d_{M,-j}(\frac{x}{|x|},y)$ we look for $b_{M,-j}$ such that
$$c_0^2(D_y^2+\lambda^2|\omb|^2/c_+^2)b_{M,-j}-\lambda^2b_{M,-j}
=d_{M,j}(\frac{x}{|x|},y).$$
The boundary conditions which $b_{M,j}$ must satisfy come from matching
with the solutions in the top and bottom hemispheres in
order to get a $C^1$ function.  They are
\begin{align}\label{eq:bc1}
b_{I,-j}(x/|x|,0)e^{i\lambda \omn y_M/c_+}
+e^{-i\lambda \omn y_M/c_+}
C_{Rj}(x/|x|)(\sin s_{R_0})^{-j}/2i\lambda c_+&=b_{M,-j}(x/|x|,y_M)\nonumber\\
i\lambda \omn/c_+[b_{I,-j}(x/|x|,0)e^{i\lambda \omn y_M/c_+}
-e^{-i\lambda \omn y_M/c_+}
C_{Rj}(x/|x|)(\sin s_{R_0})^{-j}/2i\lambda c_+]
&=
b_{M,-j}'(x/|x|,y_M)\nonumber\\
e^{-i\lambda \sqrt{1/c_-^2-|\omb|^2/c_+^2}y_M}C_{Tj}(x/|x|)
(\sin s_{T_0})^{-j}/2i\lambda c_-&=b_{M,-j}(x/|x|,-y_M)\nonumber \\
i\lambda (1/c_-^2-1/c_+^2|\omb|^2)^{1/2}
e^{-i\lambda \sqrt{1/c_-^2-|\omb|^2/c_+^2}y_M}
 C_{Tj}(x/|x|)(\sin s_{T_0})^{-j}/2i\lambda c_-&=b_{M,-j}^{\prime}(x/|x|,-y_M)
\end{align}
where $b_{I,-j+1}$ is known (it is determined by integrals over portions
of geodesics of $(\sin s)^{j-1}d_j$) and $C_{Tj}$, $C_{Rj}$
are to be determined and
are independent of $y$,
 and so can be treated as constants in solving the boundary
value problem.  They can be eliminated from this set of equations, resulting in
a boundary value problem of the type considered in Lemma \ref{l:eu},
which guarantees
us a unique solution to the problem when
$1/c_-^2>|\omb|^/c_+^2$.  This then determines $b_{R,-j}$
and $b_{T,-j}$, since $C_{T,j}$ and $C_{R,j}$ are determined by $b_{M,-j}$.

In order to ensure that our function will be $C^1$ at $y=y_M$ and at $y=-y_M$
we will add an additional term $\tilde{b}_{M,j}$ whose
total contribution will be of order
$|z|^{-j-1}$.
  Let
$\chi \in C_c^{\infty}(\Real)$, $\chi(t)=1$ for $|t|<1$ and $\chi(t)=0$ for
$|t|>2$.  Let
\begin{align*}
\beta_{Uj}&=\lim_{y \downarrow y_M}
\left(
e^{i\lambda y \omn/c_+}(b_{I,-j}(z/|z|)
+e^{-i\lambda y \omn/c_+}b_{R,-j}(z/|z|)\right)
-b_{M,-j}(x/|x|,y_M)\\
\gamma_{Uj}&=
\lim_{y \downarrow y_M}
\left(
\frac{\partial}{\partial y}(e^{i\lambda y \omn/c_+}(b_{I,-j}(z/|z|))
+\frac{\partial}{\partial
y}(e^{-i\lambda y \omn/c_+}b_{R,-j}(z/|z|))\right)
-b_{M,-j}'(x/|x|,y_M)\\
\beta_{Lj}&= \lim _{y\uparrow -y_M}\left(e^{i\lambda y
(1/c_-^2-|\omb|^2/c_+^2)^{1/2}}
b_{T,-j}(z/|z|)\right)-b_{M,-j}(x/|x|,-y_M)\\
\gamma_{Lj}&=\lim _{y\uparrow -y_M}\frac{\partial}{\partial y}\left(
e^{i\lambda y
(1/c_-^2-|\omb|^2/c_+^2)^{1/2}}
b_{T,-j}(z/|z|)\right)-b_{M,-j}^{\prime}(x/|x|,-y_M)
\end{align*}
Note that by our choice of $b_j$, $\beta_{Uj}$,
$\gamma_{Uj}$, $\beta_{Lj}$ and $\gamma_{Lj}$ all have leading order
$|z|^{-1}$.  Now, let
\begin{equation}
\tilde{b}_{M,-j}= \left(
\chi \left(\frac{3(y-y_M)}{y_M}\right)[\beta_{Uj}+(y-y_M)\gamma_{Uj}]
+\chi \left(\frac{3(y-y_M)}{y_M}\right)[\beta_{Lj}+y\gamma_{Lj}]\right)|z|.
\end{equation}
For $|y|<y_M$, this determines the approximate solution of the form
(\ref{eq:mid}).

If $1/c_-^2-|\omb|^2/c_+^2<0$, then we use a slightly different method for
finding the approximate solution when $y\leq y_M$.  Here, in
a manner similar to that used for $|y|<y_M$ above, we solve away the
error term by using an approximation of the form
$$\sum e^{i\lambda x\cdot \omb/c_+}(b_{L,-j}(\frac{x}{|x|},y)
+ \tilde{b}_{L,-j}(\frac{x}{|x|},y)|z|^{-1})|z|^{-j}$$
where
$$c_0^2(\lambda^2|\omb|^2/c_+^2+D_y^2)b_{L,-j}-\lambda^2b_{L,-j}
=d_{-j},$$
$d_{-j}$ is the coefficient of $|z|^{-j}e^{i\lambda x\cdot \omb/c_+}$
in the error term
and is
exponentially decreasing
in $y$ when $y<0$, and $b_{L,-j}$ is square integrable on $(-\infty,y_M]$.
We need in addition a boundary term at $y=y_M$, and this is
provided by the first two equations of (\ref{eq:bc1}).  An argument
like that of Lemma \ref{l:eu} shows that there is a unique solution to this
problem.
As in the previous case, $\tilde{b}_{L,-j}$ is chosen to improve the
regularity at $y=y_M$.

\begin{remark} We remark that this construction can
be carried out, with some minor modifications, for sound speeds $c_1=c+d$,
where $c-c_0$ has an asymptotic expansion of the type (\ref{eq:catinfinity}),
and $d$ is supported in $|y|<y_M$, with
$d\sim \sum_{j\geq 0}|z|^{-j}d_j(x/|x|,y)$.
\end{remark}

\subsection{Approximate Poisson Operator Near its Singularities}\label{ss:poas}

For $\omn>0$,
it remains to describe the approximation of the (partial) Poisson
operator near $z/|z|=\theta=(-\om, \omn)$ and,
if $1-c_-^2|\omb|^2/c_+^2>0$, near $z/|z|=\theta =(-c_-\omb/c_+,
-\sqrt{1-c_-^2|\omb|^2/c_+^2})$.  The
approximation in these regions contributes to the scattering matrix.
As these
two are quite similar, we will concentrate on the first, using the techniques
of \cite[Section 3]{smagailess}.

Let $w=(w_1,...,w_n)=(w',w_n)\in \Real^n$, and rotate the coordinate
system so that $\omega$ is the north pole.  Denote by
$I_{\tau}^{\gamma, \alpha}$ the class of operators
whose Schwartz kernel
 can, near the south pole, be written as a Schwartz function
plus a term of the form
\begin{equation}
\int_0^{\infty}\int \left( \frac{1}{S|w|}\right)^{\gamma}
S^{\alpha}e^{i\tau(Sw'\cdot \mu-\sqrt{1+S^2}|w|)}
a\left(\frac{1}{S|w|},S,\mu\right) dSd\mu
\end{equation}
with $a\in C_c^{\infty}([0,\epsilon)\times [0,\epsilon)\times
\Sphere^{n-2})$.  From the results of \cite{smagailess}, this class
is asymptotically complete in $\gamma$.  Moreover, a stationary phase
computation which
can be found in \cite{smagailess}
 shows that away from the south pole this is equivalent
to the class of operators whose kernel is of the form
$e^{i\tau w_n}b$, with $b$ a polyhomogeneous symbol in
$|w|$ (and smooth in $w/|w|$) of order $-\gamma+(n-1)/2$.  We recall
below
some additional facts about the operators $I^{\gamma,\alpha}_{\tau}$ from
\cite{smagailess}.

We recall from \cite{smagailess}
\begin{prop}\label{prop:sdo}
If $u(w,\om)\in I_{\tau}^{\gamma, \alpha}$ and $f\in C^{\infty}(\Sphere^{n-1}
\times \Sphere^{n-1})$, then
$$e^{i\tau|w|}\int u(|w|\theta,\om)f(\theta,\om)d\theta d\om$$
is a smooth symbolic function in $|w|$ of order $-1-\alpha$ and its lead
coefficient is $|w|^{-1-\alpha}\langle K ,f\rangle$, where $K$ is the
pull-back of the Schwartz kernel of a pseudo-differential operator
of order $\alpha-\gamma-(n-2)$ by the map $\theta\mapsto -\theta$.
The principal symbol of $K$ determines and is determined by
the lead term of the symbol, $a(t,S,\mu)$ of $u$ as $S\rightarrow 0+.$
\end{prop}

>From  \cite[Propositions 3.1 and 3.2]{smagailess}, we have
\begin{prop}
If $u\in I^{\gamma,\alpha}_{\tau}$ and $V\sim \sum_{j\geq 2}|w|^{-j}v_{j}
(w/|w|)$, then
$$(\Delta-\tau^2)u\in I^{\gamma+1,\alpha+1}$$
and $$Vu\in I^{\gamma+2,\alpha+2}.$$
\end{prop}

>From Lemma 3.2 of \cite{smagailess},
\begin{lemma}
If $u\in I^{\gamma,(n-3)/2}_{\tau}$, $V\sim \sum_{j\geq 2}|w|^{-j}v_{j}
(w/|w|)$, and $(\Delta+V-\tau^2)u\in I^{\gamma+1,\frac{n-1}{2}}_{\tau}$,
then
$(\Delta+V-\tau^2)u\in I^{\gamma+1,\frac{n+1}{2}}_{\tau}$.
\end{lemma}
Again from \cite{smagailess}
\begin{prop}
If $u\in I^{\infty,\alpha}_{\tau}=\cap_{\gamma}I^{\gamma,\alpha}_{\tau}$,
then $u=e^{-i\tau|w|}f(w)$, with $f$ a classical symbol of order $-\alpha
-1$.
\end{prop}

We proceed as described in \cite[Section 3]{smagailess}.  Using the
first part of the construction, we have an approximation of the Poisson
operator that blows up as $\theta=z/|z|$ approaches $(-\omb, \omn)$
or $(-c_- \omb /c_+, -\sqrt{1-c_-^2|\omb|^2/c_+^2})$.
Recalling that $\omn>0$,
near $\theta =(-\omb, \omn)$, we find an approximation of the Poisson
operator of the form
 $u\cal{R}_+$, where $u\in I^{-(n-1)/2,(n-3)/2}_{\lambda/c_+}$,
and $\cal{R}_+:(\omb,\omn)\mapsto (\omb,-\omn)$.
For $\theta$ near $(-c_-\omb/c_+,-\sqrt{1-c_-^2|\omb|^2/c_+^2})$,
the approximation of the Poisson operator is of the form $u\cal{R}_-$,
where $u\in I^{-(n-1)/2,(n-3)/2}_{\lambda/c_-}$ and $\cal{R}_-:
(\omb,\omn)\mapsto (c_-\omb /c_+,\sqrt{1-c_-^2|\omb|^2/c_+^2})$.
Putting all of this together, we get an approximation $\tilde{P}_0$
to the (partial)
Poisson operator with a remainder term $(c^2\Delta -\lambda^2)\tilde{P}_0$
that is in $\langle z \rangle ^{-\infty}
L^2$, and is Schwartz after
multiplication by a function $\phi \in C_b^{\infty}(\Real)$ which vanishes
for $|y|<y_M+1$.  This error is then solved away by applying
$(\Delta-(\lambda-i0)^2c^{-2})^{-1}c^{-2}$, and it is for this that we need the
results of Section \ref{s:bai}.  Taken together, the approximation of the
(partial) Poisson operator described in this section, together with
the results of Proposition \ref{prop:sdo} and Theorem \ref{thm:smc}, show
Theorem \ref{thm:sosm}.

\subsection{Approximation of $P(\lambda)\Pi_j$, $1\leq j\leq T(\lambda)$}
\label{s:pj}
For completeness, we briefly outline how to construct an approximation
$\tilde{P}_j(\lambda)$ to $P(\lambda)\Pi_j$, $1\leq j\leq T(\lambda)$.
The approximation will have the properties
$$(c^2\Delta -\lambda^2)\tilde{P}_j \in \lzr^{-\infty}L^2(\Real^n),\;
(1-\phi(y))(c^2\Delta -\lambda^2)\tilde{P}_j \in {\cal S}(\Real^n)$$
if $\phi \in C_c^{\infty}(\Real)$, $\phi(y)=1$ when $|y|<y_M+1$,
and, distributionally as $|z|\rightarrow \infty,$
$$\tilde{P}_j(z,\lambda, \omb)=|x|^{-(n-2)/2}\left(
e^{i\kappa_j(\lambda)|x|}\delta_{\omb}(\frac{x}{|x|})f_j(y)+
e^{-i\kappa_j(\lambda)|x|}h(\omb,\frac{x}{|x|})f_j(y)\right)
+{\cal O}(|x|^{-n/2}\lyr^{-\infty}).$$
Here $\omb\in \Sphere^{n-2}$.

For the construction of $\tilde{P}_j$, we begin with
$\Phi_j=e^{i\kappa_j(\lambda)x\cdot \omb}f_j(y)$, which is, up to a
constant multiple depending on $n$ and $\kappa_j$, the ($j$th partial)
Poisson operator for $c_0^2\Delta$.  We have
$$(c^2\Delta -\lambda^2)\Phi_j \sim \sum_{k\geq 2}|x|^{-k}e^{i\kappa_j
(\lambda)x\cdot \omb}d_{-k}(x/|x|,y),$$
where $d_{-k}$ is smooth in $x/|x|$, exponentially decreasing in $y$, and
smooth in $y$ when $|y|>y_M$.

To solve away the error with $k=2$, we write
$$d_{-2}(x/|x|,y)=d_{-2,1}(x/|x|,y) + d_{-2,2}(x/|x|)f_j(y)$$
where $\int d_{-2,1}(x/|x|,y)f_j(y)c_0^{-2}(y)dy=0.$  Since
$d_{-2,1}$ is orthogonal to $f_j$, we can find $g_{-2}$ such that
$$(c_0^2(D_y^2+\kappa_j^2)-\lambda^2)g_{-2}(x/|x|,y)=d_{-2,1}(x/|x|,y)$$
and use $e^{i\kappa_j(\lambda)x\cdot \omb}|x|^{-2}g_{-2}(x/|x|,y)$ to
solve away the $e^{i\kappa_j(\lambda)x\cdot \omb}|x|^{-2}d_{-2,1}(x/|x|,y)$
error, up to a term vanishing one order faster at infinity.  We note
that $g_{-2}$ is exponentially decreasing in $y$, since $d_{-2,1}$ is,
and that this term does not contribute anything to the scattering
matrix, since distributionally it is ${\cal O}(|x|^{-(n+2)/2})$.

To solve away the error
$|x|^{-2}e^{i\kappa_j(\lambda)x\cdot \omb}d_{-2,2}(x/|x|)
f_j(y)$, we use the techniques of \cite{smagailess} in the $x$ variables only.
That is, essentially as in our construction of $\tilde{P}_0$, we
solve transport equations along geodesics on $\Sphere^{n-2}$ beginning at
$x/|x|=\omb$.  Near $x/|x|=-\omb$, we must use the second ansatz, as in
Section \ref{ss:poas} or \cite{smagailess}.

The subsequent errors are solved away in exactly the same manner, resulting
in an approximation as claimed.

\section{End of proof of structure results}\label{s:bai}

In order to finish the proof of Theorem \ref{thm:sosm},
we must tie up a number of
loose ends.  These include showing the existence of a function $u$
as in the definition of the Poisson operator (Definition \ref{d:po}),
proving Proposition \ref{p:npoe}, and showing that the
construction of the (partial) Poisson operators in the
previous sections captures all the singularities of the scattering matrix.  As
these use similar techniques, we give the proofs in this section.

For $\lambda \in \Real \setminus \{0\}$, $\ep>0$, the limit
\begin{multline*}
\lim_{\delta \downarrow 0} \lzr^{-1/2-\ep}(c^2\Delta
-(\lambda-i\delta)^2)^{-1}c^2
\lzr^{-1/2-\ep}
=  \lim_{\delta \downarrow 0} \lzr^{-1/2-\ep}
(\Delta -c^{-2}(\lambda-i\delta)^2)^{-1}
\lzr^{-1/2-\ep}\\
=  \lzr^{-1/2-\ep}
(\Delta -c^{-2}(\lambda-i0)^2)^{-1}
\lzr^{-1/2-\ep}
\end{multline*}
as an operator on $L^2(\Real^n)$ exists in the norm topology (\cite{bdm-m,
db-pi}).  In this section, we study further properties of
$(\Delta -(\lambda-i0)^2c^{-2})^{-1}$ when it is applied to a function
$f\in \lzr^{-\infty}L^2(\Real^n)$, $(1-\phi(y))f\in \sch (\Real^n)$ for
some $\phi\in C_c^{\infty}(\Real)$.  In particular, we are interested in
the asymptotics at infinity of the resulting function.

For simplicity of exposition, we shall assume $\lambda>0$ throughout
 this section.  The results for $\lambda<0$ can be proved in a similar way.

\subsection{``Outgoing'' solutions}
 The main result of this section is Proposition \ref{prop:rio}.  This
proposition, combined with the approximation of the Poisson operator in the
previous section, shows the existence of a function $u$ with the properties
given in Definition \ref{d:po}.
Taken together with Proposition \ref{p:powd}, this shows the existence
and uniqueness of the Poisson operator.

In proving this proposition, as well as in many others, we shall use the
fact that if
$$u=\lim_{\ep \downarrow 0}(\Delta -c^{-2}(\lambda-i\ep)^2)^{-1}f=
(\Delta -c^{-2}(\lambda-i0)^2)^{-1} f,$$ then
\begin{equation}\label{eq:rr}
u=(\Delta -(\lambda-i0)^2 c^{-2}_0)^{-1}(Vu+f)
\end{equation} where
\begin{equation}\label{eq:v}
V=\lambda^2(c^{-2}-c_0^{-2}).
\end{equation}
Additionally,
\begin{equation}\label{eq:rc0}
(\Delta -(\lambda-i0)^2c_0^{-2})^{-1}g(z)= (2\pi)^{1-n}
\int e^{ix\cdot \xi}
((D_y^2+|\xi|^2-(\lambda-i0)^2c_0^{-2})^{-1}\hat{g}(\xi,\cdot))(y)d\xi
\end{equation}
where $\hat{g}(\xi,y)=\int e^{-ix\cdot \xi}g(x,y)dx$ is the Fourier
transform in the $x$ variables only.  We shall repeatedly use this notation
for the Fourier transform in the $x$ variables only, and $V$ is given by
(\ref{eq:v}).

We make several remarks about the operator
$(D_y^2+t^2-(\lambda-i0)^2c_0^{-2})^{-1}$.  As an operator from
$(D_y^2+t^2-(\lambda-i0)^2c_0^{-2})^{-1}:\lyr^{-1/2-\ep}L^2(\Real)\rightarrow
\lyr^{1/2+\ep}L^2(\Real)$ for any $\ep >0$ it is smooth for $|t|<\lambda/c_+$,
$t$ away from $\lambda/c_{\pm}$.  Near $t=\lambda/c_+,$
$(\lambda^2/c_+^2-t^2)^{1/2}(D_y^2+t^2-(\lambda-i0)^2c_0^{-2})^{-1}$ is a
smooth
function of $(\lambda^2/c_+^2-t^2)^{1/2}$, and, if $c_->c_+$, near
$t=\lambda/c_-$ it is a smooth function of $(\lambda^2/c_-^2-t^2)^{1/2}$.

We shall use the following lemma in the proof of Proposition
\ref{prop:rio}.
\begin{lemma}\label{l:tsm}
If
$u = \lim_{\epsilon \downarrow 0}
(\Delta-(\lambda-i\epsilon)^2c^{-2})^{-1}f$ and
 $f\in \langle z \rangle ^{-\infty}L^2$, then
$\widehat{Vu}(\xi,y)$ and all its derivatives with respect to $\xi/|\xi|$
are, for
$\epsilon>0$, in
$H^{J-1/2-\epsilon-\beta}(\Real^{n-1}_{\xi};
\langle y \rangle ^{-\beta}L^2(\Real_y))$
for $0\leq \beta \leq J-1/2-\epsilon$.
\end{lemma}
\begin{pf}
Throughout the proof $\epsilon>0$ is small, and may change from
line to line.

We have
$(\Delta -\lambda^2/c_0^2)u=Vu+f$.
Then
\begin{equation}\label{eq:uid}
 \widehat{Vu}(\xi,y)=C
\int e^{-i x \cdot \xi}V(x,y)\int e^{ix\cdot \eta}
((D^2_y+|\eta|^2-c^2_0(\lambda-i0)^2)^{-1}[\hat{f}(\eta,\cdot)+
\widehat{Vu}(\eta,\cdot)])(y)d\eta dx.
\end{equation}
Let $D_{(\xi/|\xi|)_l}$ stand for a derivative tangent to $|\xi|=$constant.
Then
\begin{multline}\label{eq:com1}
\left(  D_{(\xi/|\xi|)_l}\right)^j\int
e^{-i x \cdot \xi}V(x,y)\int e^{ix\cdot \eta}
((D^2_y+|\eta|^2-c^{-2}_0(\lambda-i0)^2)^{-1}
\widehat{Vu}(\eta,\cdot))(y)d\eta dx \\
=\int
e^{-i x \cdot \xi}V(x,y)\int e^{ix\cdot \eta}
((D^2_y+|\eta|^2-c^{-2}_0(\lambda-i0)^2)^{-1}
\left(  D_{(\eta/|\eta|)_l}\right)^j
\widehat{Vu}(\eta,\cdot))(y)d\eta dx   \\
+\sum_{|\alpha|<j} \int
e^{-i x \cdot \xi}V_{\alpha}(x,y,\xi)\int e^{ix\cdot \eta}
((D^2_y+|\eta|^2-c^{-2}_0(\lambda-i0)^2)^{-1}\left(  D_{(\eta/|\eta|)}\right)^
{\alpha}
\widehat{Vu}(\eta,\cdot))(y)d\eta dx.
\end{multline}
Here $D_{(\eta/|\eta|)}^{\alpha}$ is a derivative of order $|\alpha|$ in
the $\xi/|\xi|$ variables and
 \begin{equation}
|D^{\beta}_x D^k_yD^{\gamma}_{\xi} V_{\alpha}(x,y,\xi)|
\leq C(1+|z|)^{-J-|\beta|-k}
\end{equation}
for large $|z|$, where if $k>0$, we need also require $|y|>y_M$.  That
is, $V_{\alpha}$ is a symbol in $x$, a property we shall use.
 Let $\dxt$  be the Laplacian in the $\xi/|\xi|$ variables.  If $0<t<1/2$,
then
\begin{multline}\label{eq:com2}
\left(  \dxt \right)^{t}\int
e^{-i x \cdot \xi}V(x,y)\int e^{ix\cdot \eta}
((D^2_y+|\eta|^2-c^{-2}_0\lambda^2)^{-1}
\widehat{Vu}(\eta,\cdot))(y)d\eta dx \\
=\int
e^{-i x \cdot \xi}V(x,y)\int e^{ix\cdot \eta}
((D^2_y+|\eta|^2-c^{-2}_0\lambda^2)^{-1}\left(  \Delta_{\eta/|\eta|}\right)^t
\widehat{Vu}(\eta,\cdot))(y)d\eta dx  \\
+\int
e^{-i x \cdot \xi}V_{t}(x,y,\xi)\int e^{ix\cdot \eta}
((D^2_y+|\eta|^2-c^{-2}_0\lambda^2)^{-1}
\widehat{Vu}(\eta,\cdot))(y)d\eta dx
\end{multline}
where $V_{t}(x,y,\xi) $ satisfies, for large $|z|$,
$$|D_{x}^{\alpha}D_y^{m}D_{\xi}^{\gamma}V_t(x,y,\xi)|
\leq C_{\alpha, t,\gamma, m}(1+|z|)^{-J-1+2t-|\alpha|-m}$$
where we require $m=0$ if $|y|\leq y_M$.

Because of the decay properties of $V$, $\widehat{Vu}(\eta,y)
\in H^{J-1/2-\epsilon
-\beta}(\Real^{n-1}_{\eta};\langle y \rangle ^{-\beta}L^2(\Real_y))$
 for $J-\epsilon -1/2\geq
\beta \geq 0$.  Of course, $\hat{f}(\eta,y) \in H^{\infty}
(\Real^{n-1}_{\eta};\langle y \rangle ^{-\infty}L^2_y)$.

Suppose for some $t>0$, $( \dxt )^{t/2}\widehat{Vu}
\in H^{J-1/2-\epsilon -\beta}(\Real^{n-1}_{\xi};\langle y \rangle ^{-\beta} L^2
(\Real_y))$
for $J-\epsilon -1/2\geq
\beta \geq 0$.  Then using equations (\ref{eq:uid}),
(\ref{eq:com1}) and (\ref{eq:com2})
we obtain  $( \dxt)^{(t+2J-2-2\epsilon)/2}\widehat{Vu}
\in H^{J-1/2-\epsilon -\beta}(\Real^{n-1}_{\xi};\langle y \rangle ^{-\beta}
L^2(\Real_y))$.
Since $J>1$, the proof follows by induction.
\end{pf}

\begin{lemma}\label{l:spdo}
If $w,\eta\in \Real^m$, and
$$|D^{\alpha}_{z}D_{\eta}^{\beta}a(w,\eta)|\leq
C_{\alpha\beta}(1+|w|)^{-|\beta|}$$ for all multiindices $\alpha, \beta$, then
the operator $A$ defined
by
$$(Ag)(\eta)=\int e^{iw\cdot \eta} a(w,\eta)g(w)dw$$
is a continuous map $$A:\langle w \rangle^{\gamma}L^2(\Real^m)\rightarrow
H^{-\gamma}
(\Real^m).$$
\end{lemma}
\begin{pf}  We can write
$$(Ag)(\eta)=\int e^{iw\cdot (\eta-\eta')}a(w,\eta)\check{g}(\eta')d\eta'dw$$
with $\check{g}\in H^{-\gamma}(\Real^m)$
 the inverse Fourier transform on $\Real^m$.  Since $a$ is
a symbol of order $0$ in $w$, this is a pseudodifferential operator acting
on $\check{g}$ and the result follows from standard pseudodifferential operator
theory.
\end{pf}

\begin{lemma}\label{l:hgzg}
Let $a(w,\eta)$ be a symbol of order $0$ in $w=(\overline{w},w_{m+1})
\in \Real^{m+1}$, and let
 $h(\eta)
\in H^{-\gamma}
(\Real^m)$ be supported away from $t^2=|\eta|^2$, where $t$ is a nonzero
constant.  Then
$$(Ah)(w)=\int e^{i\overline{w}\cdot \eta +iw_{m+1}\sqrt{t^2-|\eta|^2}}
 a(w,\eta)
h(\eta) d\eta \in \langle w \rangle ^{\gamma + 1/2}L^2(\Real^{m+1})$$
provided $\gamma>0$.
\end{lemma}
\begin{pf}
Let $g\in \langle w \rangle ^{-\gamma -1/2}L^2(\Real^{m+1})$.  Then
\begin{equation}\label{eq:gp}
(g,Ah)=(\int \overline{a}(w,\eta)e^{-i\overline{w}\cdot \eta
-iw_{m+1}\tau}g(w)dw_{\rest \tau=\sqrt{t^2-|\eta|^2}},h(\eta)).
\end{equation}
By the previous lemma and the restriction properties of elements
of Sobolev spaces,
$$\int \overline{a}(w,\eta)e^{-i\overline{w}\cdot \eta
-iw_{m+1}\tau}g(w)dw_{\rest
\tau=\sqrt{t^2-|\eta|^2}}\in H^{\gamma}$$ if $\gamma>0$.  The pairing
(\ref{eq:gp}) is
then well-defined for all such $g$, and
$Ah\in \langle w \rangle^{\gamma +1/2}L^2(\Real^{m+1})$.
\end{pf}

\begin{prop}\label{prop:rio}
If $f\in \lzr^{-\infty}L^2(\Real^n)$ and $J\geq 2$, then
$u=(\Delta-(\lambda-i0)^2c^{-2})^{-1}f$ is outgoing in the sense of Definition
\ref{d:outgoing}.
\end{prop}
\begin{pf}

We use (\ref{eq:rr}) and (\ref{eq:rc0}).  Choose $\Psi \in C_c^{\infty}(\Real)$
to be $1$ for  $|\xi | \leq \lambda^2/c_+^2$ and supported in
a slightly larger neighborhood.
We use the fact that we can write, for $\delta >0$, $(1-\Psi(\xi))
((D_y^2+|\xi|^2-(\lambda-i\delta)^2c_0^{-2})^{-1}$ as a sum of an operator
bounded on $L^2(\Real)$ and an operator involving
projection onto the discrete spectrum.  The part corresponding to
an operator bounded on $L^2$ gives an element of $L^2$.
For the part corresponding to
the discrete spectrum, near the poles
we move the contour of integration in $|\xi| $ slightly into the upper
half-plane when $x\cdot \xi \geq 0$ and into the lower half-plane for
$x\cdot\xi \leq 0$, obtaining
\begin{equation*}
\int e^{ix\cdot \xi}(1-\Psi(|\xi|))
((D_y^2+|\xi|^2-(\lambda-i0)^2c_0^{-2})^{-1}(\widehat{Vu}+\widehat{f})
(\xi,\cdot))(y)d\xi
= \sum_1^{T(\lambda)} e^{-i\kappa_j(\lambda)|x|} b_j(\frac{x}{|x|})f_j(y) +
\tilde{u}_1
\end{equation*}
where $\tilde{u_1}, \frac{\partial}{\partial |z|} \tilde{u}_1
\in L^2(\Real^n)$.  We remark here that $b_j\in C^{\infty}(\Sphere^{n-1})$,
thanks to the fact that it originates from a stationary phase applied to
functions smooth in $\xi/|\xi|$ (Lemma \ref{l:tsm}) at
$|\xi|=\kappa_j(\lambda)$.  We refer the reader to the proof of Theorem
4.1 of \cite{tsm} for greater detail in a very similar computation.

When $|y|<C$ for some constant $C$, we have that away from $|\xi|=\lambda/c_+$,
$$w(\xi,y)=
\Psi(|\xi|)((D_y^2+|\xi|^2-(\lambda-i0)^2)^{-1}(\widehat{Vu}+\widehat{f})
(\xi,\cdot))(y)\in L^2(\Real^{n-1}_{\xi})$$
with a norm independent of $y$, $|y|<C$.
Near $|\xi|=\lambda/c_+$, $w(\xi,y)\in L^p(\Real^{n-1})$ for any
$p<2$, again with norm
independent of $y$, $|y|<C$, so using the mapping properties of the
Fourier transform we have
$$\int e^{ix\cdot \xi}w(\xi,y)d\xi_{\rest \{|y|<C\}}\in \langle x \rangle
^{\ep} L^2(\Real^{n-1}_x\times [-C,C]_y).$$

When $y>y_M$, for $|\xi|\leq \lambda/c_+$ we have
$$
((D_y^2+|\xi|^2-(\lambda-i0)^2c_0^{-2})^{-1}\widehat{Vu}+\widehat{c^2f}
(\xi,\cdot))(y) = h_1(y)e^{-i\sqrt{\lambda^2/c_+^2-|\xi|^2}}+h_2(y,\xi).$$
Since
$h_2(y,\xi)\in L^2(\Real^{n-1}_{\xi};\lyr^{-J+3/2-\ep}L^2(\Real_y))$,
taking its inverse Fourier transform in $x$ gives an element of
$\lyr^{-J+3/2-\ep}L^2(\Real^n\cap \{y>y_M\})$, and similarly
if we take a radial derivative.

For the term
$h_1(y)e^{-i\sqrt{\lambda^2/c_+^2-|\xi|^2}}$, note that away from
$|\xi|=\lambda/c_{\pm}$, $h_1\in H^{J-1-\ep}(\Real^{n-1})$.
 Notice that
\begin{equation}\label{eq:ybig}
(\frac{\partial}{\partial |z|}+i\lambda/c_+)e^{ix\cdot
\xi-iy\sqrt{\lambda^2/c_+^2-|\xi|^2}}= (\frac{x}{|z|}\cdot \xi -
i \frac{y}{|z|}\sqrt{\lambda^2/c_+^2-|\xi|^2}+i\lambda/c_+)
e^{ix\cdot \xi-iy\sqrt{\lambda^2/c_+^2-|\xi|^2}},
\end{equation}
that is, it vanishes to first order on the critical points of the phase.
Therefore, smoothly restricting our region of integration
to  $|\xi|<\lambda/c_+-\delta$ and away from $|\xi|=\lambda/c_-$,
an integration by parts argument shows us that we get an element of
$\lzr^{\ep}L^2(\Real^n\cap \{y>y_M\})$ for any $\ep>0$.

 Near
$|\xi|=\lambda/c_+$,
$|\xi|\leq \lambda/c_+$,
 $h_1= (\lambda^2/c_+^2-|\xi|^2)^{-1/2}b(\sqrt{\lambda^2/c_+^2-|\xi|^2})h_+
(\xi)$, where $b\in C^{\infty}$ and
$h_+(\xi)=h_{++}(\xi, \sqrt{\lambda^2/c_+^2-|\xi|^2})$,
$h_{++}\in
 H^{J-1/2-\ep}(\Real^n)$. Using a partition of
unity, we can work on coordinate patches on
which $\xi_i \not = 0$ for some $i$.  For example, let $\overline{\xi}=
(\xi_2,\xi_2,...\xi_{n-1})$ and
for $\xi_1>\delta>0$, we use
$(\overline{\xi},t=\sqrt{\lambda^2/c_+^2-|\xi|^2})$
 as our variables of
integration.
Then we need to compute
$$(\frac{\partial}{\partial |z|}+i\lambda/c_+)
\int _{t\geq 0}
e^{i(x_1\sqrt{\lambda^2/c_+^2-t^2-|\overline{\xi}|^2}+\overline{x}
\cdot
\overline{\xi}-yt)}\chi_1(\overline{\xi},t)b(t)h_{++}(\sqrt{\lambda^2/c_+^2-t^2-|\overline{\xi}|^2},
\overline{\xi},t)dt.$$   Here $\chi_1$ is a smooth
cut-off function and $\overline{x}=(x_2,x_3,...,x_{n-1})$.
 The resulting integrand
vanishes at the stationary points of the phase, so that we may integrate
by parts.  Applying Lemma \ref{l:hgzg} and a related result for the boundary
term,
we obtain an element of $\lzr^{\ep}L^2(\Real^n)$.  We can repeat this
procedure, using the partition of unity, to cover the integration
over the region with $|\xi|\leq \lambda/c_+$,
$|\xi| $ near $\lambda/c_+$.

For $|\xi|\geq \lambda/c_+$, we have
\begin{multline}\label{eq:xi>l/c_+}
\int_{|\xi|\geq \lambda/c_+}\Psi(|\xi|)e^{ix\cdot \xi}
(|\xi|^2-\lambda^2/c_+^2)^{-1/2}\left(
e^{-y\sqrt{|\xi|^2-\lambda^2/c_+^2}}
\right.\\ \left. \int_{-\infty}^y a(\sqrt{|\xi|^2-
\lambda^2/c_+^2})(\widehat{Vu}+\hat{f})(\xi,y')
\phi(y',\sqrt{|\xi|^2-\lambda^2/c_+^2})dy'
+h_{2b}(\xi,y)\right)d\xi
\end{multline}
with $a(t)$, $\phi(y',t)$ smooth in $t$ and
$|\phi(y',t)|\leq Ce^{y't}$.  The function
$\phi(y,\sqrt{|\xi|^2-\lambda^2/c_+^2})$ is in
the null space of $D_y^2+|\xi|^2-\lambda^2/c_0^2$.  Since
$h_{2b},D_yh_{2b}\in H^{J-1-\ep-\gamma}
(\Real^{n-1};\lyr^{-\gamma}L^2(\Real))$ for $J-1-\ep-\gamma\geq 0$,
its contribution is
in $\lzr^{\ep}L^2(\Real^n)$ and can be neglected.  As in the previous integral,
for the rest we
use a partition of unity to work in sets with $\xi _i \not = 0$, some $i$.
Near $\xi_1 \not =0$, we use coordinates $\overline{\xi}, t=\sqrt{|\xi|^2-
\lambda^2/c_+^2}$ almost as before.  Using the technique of the
proof of Lemma \ref{l:hgzg},
applying
$(\frac{\partial}{\partial |z|}+i\lambda/c_+)$ when $y>y_M$, integrating
by parts, and using Lemma \ref{l:spdo}, we need to check only that
\begin{equation}\label{eq:gpy}
\int\int_{t\geq 0}\int_{y_M}^{\infty}
 w_1(\overline{\xi},t,y)\int_{y_M}^ye^{-yt}\phi(y',t)w_2(\overline{\xi},t,y')dy'dt
d\overline{\xi}
\end{equation}
is finite
when $w_1\in H^{\ep/2}(\Real^{n-1};\lyr^{-1-\ep/2}L^2(\Real))$
and $w_2\in H^{J-1-2
\delta}(\Real^{n-1};\lyr^{-\delta}L^2(\Real))$ for any $\delta >0$,
where $w_2$ is supported near $t=0$ and away from
$|\overline{\xi}|=\lambda/c_+$.
Since by H\"older's inequality
$$|\int_{y_M}^ye^{-(y-y')t}w_2(\overline {\xi},t,y')dy'|
\leq C \lyr^{-\delta+1/2}(\int | w_2(\overline{\xi},t,y')|^2\langle y'
 \rangle^{2\delta}
dy')^{1/2}$$ for
any $\delta >0$, using H\"older's inequality again we see that
(\ref{eq:gpy}) is finite.

Near $|\xi|=\lambda/c_-$, $h_1=h_{-1}+h_{-2}$, where $h_{-2}\in H^{J-1-\ep}$
and can be treated by an integration by parts as before.  However,
$h_{-1}$ is the restriction of an element of $H^{J-1/2-\ep}$ to the
hypersurface $(\xi, \sqrt{\lambda^2/c_-^2-|\xi|^2})$.  Here we introduce the
variable of integration $t= \sqrt{\lambda^2/c_-^2-|\xi|^2}$; this
results in an extra factor of $t$ in the integrand and makes $t=0$ a
stationary point in of the phase.  However, the extra factor of $t$ in
the integrand is enough to allow us to integrate by parts after
applying $\frac{\partial}{\partial|z|}+i\lambda/c_+$,
and arguing as in Lemma \ref{l:hgzg}, we obtain an element of
$\lzr^{\ep}L^2(\Real^n\cap\{y>y_M\})$.

A similar technique shows that
$(\frac{\partial}{\partial
|z|}+i\lambda/c_-)(\Delta-(\lambda-i0)^2c^{-2})^{-1}f_{\rest y<-y_M}\in
\lzr^{\ep}L^2(\Real^n\cap\{y<-y_M\})$ for any
$\ep>0$.
\end{pf}

For future reference, we remark that the proof above has shown
\begin{cor}\label{c:evp}
If $f \in \lzr^{-\infty}(\Real^n)$,
then for $y\in K$, $K\subset \Real$
compact, $|x|>1$,
$$(\Delta -c^{-2}(\lambda-i0)^2)^{-1}f(z)_{\rest|x|>1, y\in K}
 =|x|^{-(n-2)/2}\sum_1^{T(\lambda)}
e^{-i\kappa_j(\lambda)|x|}b_j(x/|x|)f_j(y)+u_1$$
where $b_j \in C^{\infty}(\Sphere^{n-2})$ and $u_1\in \langle x \rangle^{\ep}
L^2(\Real^{n-1}\times K)$.
\end{cor}

\subsection{Smoothness of leading order coefficient}\label{s:smooththm}
To finish the proof of Theorem \ref{thm:sosm}, we need to prove Theorem
\ref{thm:smooth}, which we recall:
\begin{thmb}
 Let $c$ and $c_0$ satisfy the hypotheses of
Section \ref{s:aan} and hypothesis (H1) or (H2).
For any $\chi \in C_c^{\infty}(\Spc^{n-1})$, $f\in \lzr^{-\infty}L^2
(\Real^n)$,
$(1-\phi(y))f\in {\cal S}(\Real^n)$ for some $\phi \in C_c^{\infty}(\Real)$,
we have
$$\chi(z/|z|)(\Delta -c^{-2}(\lambda-i0)^2)^{-1}f
= e^{-i\lambda|z|/c}|z|^{-(n-1)/2}a_0(z/|z|) +u_1$$
where $a_0 \in C^{\infty}(\Spc^{n-1})$ and $u_1 \in \lzr^{\ep}L^2(\Real^n)$
for all $\ep >0$.
\end{thmb}
The
proof of this theorem will occupy the next two subsections.  Before
we prove it, however, we give some applications.

\vspace{3mm}
\noindent
{\em Proof of Proposition \ref{p:npoe}}.  We have
$P(\lambda)\Pi_j=\tilde{P}_j(\lambda)-(\Delta -c^{-2}(\lambda-i0)^2)^{-1}
(
\Delta -\lambda^2 c^{-2})\tilde{P}_j(\lambda)$, where $\tilde{P}_j$ are the
approximate
(partial)
Poisson operators constructed in Section \ref{s:cpo}.  Proposition
\ref{p:npoe} now follows from our construction of $\tilde{P}_j$, Proposition
\ref{prop:rio},
Corollary \ref{c:evp} and Theorem
\ref{thm:smooth}.
\endpf
\vspace{3mm}

Moreover,
if we take $f=(\Delta-\lambda^2c^{-2})\tilde{P}_0(\lambda, z,\om)$, where
$\tilde{P}_0$ is the approximation of the (partial)
Poisson operator we have constructed, then this theorem
 shows that the the only singularities of the main part of the
scattering matrix come from the structure of the approximation
of the Poisson operator we have constructed.  This then proves
Theorem \ref{thm:sosm}.

The following proposition implies Theorem \ref{thm:smooth} in one case.
\begin{prop}
Suppose $c$ and $c_0$ satisfy the hypothesis (H1).  Let $\chi \in C_c^{\infty}
(\Spc^{n-1})$ and let $f\in \lzr^{-\infty}L^2$,
$(1-\phi(y))f\in {\cal S}(\Real^n)$ for some $\phi \in C_c^{\infty}(\Real).$
Then
$$\chi(z/|z|)(\Delta -c^{-2}(\lambda -i0)^2)^{-1}f
= |z|^{-(n-1)/2}e^{-i\lambda |z|/c_0(y)}\sum_0^N |z|^{-j}a_j(z/|z|) + {\cal
O}(|z|^{-(n-1)/2-N-1})$$
where $a_j \in C^{\infty}(\Spc^{n-1})$.
\end{prop}
\begin{pf} Because $c_+=c_-$, we can recast this in the form of a perturbation
of a
simple $n$-body problem and use the fact that much is known about asymptotic
expansions of the resolvent applied to Schwartz function.
First, we note that
$$(\Delta -c^{-2}(\lambda-i0)^2)^{-1}=(\Delta -\lambda^2(c^{-2}-c_+^{-2})
-c_+^{-2}(\lambda-i0)^2)^{-1}.$$
The operator $\Delta -\lambda^2(c^{-2}_0-c_+^{-2})$ is a particularly simple
example of a class of
$n$-body operators widely studied.
The operator $\Delta -\lambda^2(c^{-2}-c_+^{-2})$, while
not quite an $n$-body operator since the potential depends on all variables, is
a perturbation that has many of the same properties we desire.

The paper \cite{vasyab}, which builds on results of \cite{g-i-s,sslaes},
shows that \begin{multline*}\chi(z/|z|)(\Delta - \lambda^2(c^{-2}_0-c_+^{-2})
-c_+^{-2}(\lambda-i0)^2 )^{-1}f \\
= |z|^{-(n-1)/2}e^{-i\lambda |z|/c_0(y)}\sum_0^N |z|^{-j}b_j(z/|z|) + {\cal
O}(|z|^{-(n-1)/2-N-1})
\end{multline*}
with $b_j$ smooth.  The proof is such that the results of \cite{vasyab} hold
with $c_0$ replaced by a sound speed $c$ of the type considered here.  Roughly
speaking, this
is because \cite{g-i-s} requires that the operator
$\Delta +V_1$ (where for us $V_1=  - \lambda^2(c^{-2}_0-c_+^{-2})$ or
$V_1=  - \lambda^2(c^{-2}-c_+^{-2})$ ) satisfy a Mourre estimate and
some regularity and decay
 properties, both of which are satisfied for either $V_1$.
Vasy remarks already that the results of \cite[Section 2]{vasyab}, which
are local versions of results of \cite{sslaes}, will hold in our case.
Then the results of \cite[Section 3]{vasyab} hold for our case, since
just as in that paper we can argue that from the results of \cite{g-i-s}
and \cite{sslaes} that $$\WF^*_{\sca}\left(
\chi(z/|z|)(\Delta -c^{-2}(\lambda -i0)^2)^{-1}f
\right)
\subset R^+_{\lambda/c_+},$$ and the asymptotic expansion follows from
\cite[Proposition 2.8]{vasyab} and the remarks made there.
(Roughly speaking, $\WF_{\sca}(u)$ provides a microlocal description of
the lack of decay of $u$-- see \cite{sslaes} for a definition of
the scattering wave front set, $\WF_{\sca}$, and $R^+_{\tau}$.)
\end{pf}

Before giving the proof in the case $c_+\not = c_-$, we give
some explanation as to why the proof of \cite{vasyab} does not apply
in this case.  The argument of \cite{vasyab} uses very strongly
the fact that, for $f\in {\cal S}$, $$\WF_{\sca}^*\left(\chi(z/|z|) (\Delta
-\lambda^2(c^{-2}_0-c_+^{-2})-c_+^{-2}(\lambda -i0)^2)^{-1} f\right)
\subset R^+_{\lambda/c_+}$$
where $\chi \in C_c^{\infty}(\Spc^{n-1}).$

This is not, however, true in general when $c_+\not = c_-$.
The results of
\cite[Theorem 4.1]{tsm} show that in general for $f\in {\cal S}(\Real^n)$,
\begin{multline}
\left((\Delta -(\lambda-i0)^2c_0^{-2})^{-1}f \right)(z)
= e^{-i \lambda |z|/c_+}|z|^{-(n-1)/2}a_0(z/|z|)\\
+ |z|^{-(n+1)/2}\left(e^{-i\lambda |z|/c_+}a_1(z/|z|)+e^{-i\lambda|x|/c_-
-iy\lambda \sqrt{1/c_+^2-1/c_-^2}}b_1(z/|z|)\right) +
{\cal O}(|z|^{-(n+3)/2})
\end{multline}
when $0<\ep <y/|z|<(1-c_+^2/c_-^2)^{1/2}-\ep$ and $|z|\rightarrow \infty$.
If $b_1\not = 0$, then the scattering wave front set is not
contained in $R_{\lambda/c_+}^+$.
This means that the scattering wave front set of
$\chi (\Delta -(\lambda -i0)^2/c^2)^{-1}f $ is in general more
complicated than the $c_+=c_-$ case (and
unknown, to the best of our knowledge), and the
techniques of \cite{sslaes,vasyab} cannot be immediately applied.
Similar differences can be seen in the resolvent estimates of
\cite[Theorem 1.1]{g-i-s} for the
$n$-body problem, and \cite[Theorem 3.1]{hi}, for a particular
stratified medium.

Instead, we take a different approach.

\subsection{Proof of Theorem \ref{thm:smc} in case hypothesis (H2)
holds}

In order to prove Theorem \ref{thm:smc} when hypothesis (H2) holds, we will
make heavy use of equations (\ref{eq:rr}) and (\ref{eq:rc0}).  We use the
fact that the more rapidly $g$ decays at infinity, the more we can say
about $(\Delta -(\lambda-i0)^2c_0^{-2})^{-1}g$,
using equation (\ref{eq:rc0}).
To take advantage of this, roughly speaking, we find
approximations $w$ of $u=(\Delta -(\lambda-i0)^2c^{-2})^{-1}f$ so that
$(\Delta -\lambda^2/c_0^2)w=Vu+e$, with the error $e$ decaying faster than
$Vu$ does, and $w$ decaying faster than $u$.  Then
$u=w-(\Delta -(\lambda-i0)^2c_0^{-2})^{-2}(f+e)$ (compare (\ref{eq:rr})).
The better rate of decay of $e$ improves our knowledge of $u$.

In practice, the proof is somewhat more complicated.  We study the behaviour
at infinity of $\chi(z/|z|)u$, $\chi\in C_c^{\infty}(\Spc^{n-1})$, and
we introduce a ``microlocal'' cut-off $\Psi(D_x)$ so that
$\chi\Psi(D_x)u$ has the same leading behaviour as $\chi u$ at infinity
(Lemma \ref{l:leadorder}), but is easier to understand.

Lemma \ref{l:leadorder}
 shows that, for suitable
$\Psi$, $\chi(1-\Psi(D_x))u\in \lzr^{\ep-1/2}L^2(\Real^n)$,
so that $\chi\Psi(D_x)u$ captures the leading behaviour of $\chi u$.
Lemmas \ref{l:cp}-\ref{l:dr} are preliminaries, and
Lemmas \ref{l:ps1}-\ref{l:gft}
 are used to make successive approximations of $\chi u$
in the proof of Proposition \ref{prop:slot}, which shows that $\chi \Psi(D_x)u$
has an asymptotic expansion with smooth coefficient in the leading order term.

The next lemma shows us that, for suitable $\Psi$ and $\chi$, $\chi \Psi(D_x)u$
is the leading order term of $\chi u$.
\begin{lemma} \label{l:leadorder} Let $\chi\in C^{\infty}_c(\Spc^{n-1}),$
$\Psi\in C^{\infty}_b(\Real^{n-1})$, and suppose that if $(x/|z|,y/|z|)\in
\supp \chi$ and $y>0$, then $\pm \lambda x/(c_+|z|)\not \in \supp(1-\Psi)$,
and if $(x/|z|,y/|z|)\in
\supp \chi$ and $y<0$, then $\pm \lambda x/(c_-|z|)\not \in \supp(1-\Psi)$.
If $J\geq 4$,
$u = (\Delta -(\lambda-i0)^2c^{-2})^{-1}f$,
and $f\in \lzr^{-\infty}L^2(\Real^n),$ then
$$\chi(z/|z|)(1-\Psi(D_x))u, \;
\frac{\partial}{\partial |z|}\chi(z/|z|)(1-\Psi(D_x))u
 \in \lzr^{-1/2+\ep}L^2(\Real^n)$$
for any $\ep >0$.
\end{lemma}
\begin{pf}
We give the proof for $\chi$ supported in $y>0$, as the proof for $\chi$ with
support in $y<0$ is quite similar.

The proof closely resembles that of Proposition \ref{prop:rio}.  Let
$\Psi_1 \in C_b^{\infty}(\Real)$ be such that $\supp \Psi_1(t)
\subset \{ |t|>\lambda/c_+\}$, $\supp (1-\Psi_1)\subset \{|t|\leq \lambda/c_+
+\delta\}$, some small $\delta >0$.
Then, by the same type of arguments as in the proof of Proposition
\ref{prop:rio}, since the eigenfunctions of $D_y^2+c_0^{-2}\lambda^2$
are exponentially decreasing in $y$, we have
$$\chi(z/|z|)\int e^{ix\cdot \xi}(1-\Psi(\xi))\Psi_1(|\xi|)
(D_y^2+|\xi|-c_0^{-2}(\lambda - i0)^2)^{-1}
(\widehat{Vu}(\xi,\cdot)+\hat{f}(\xi,\cdot))(y)d\xi
\in \langle z\rangle ^{-1/2+\ep}L^2(\Real^n).$$
Moreover, for $y>y_M$, using (\ref{eq:ybig}), we see that when
$z$ is restricted to the support of $\chi(z/|z|)$ there
are no stationary points of the associated phase $x\cdot \xi -y\sqrt{\lambda^2
/c_+^2-|\xi|^2}$ on the support of $1-\Psi(\xi).$  Therefore, when the
integrand is supported in
$|\xi|<\lambda/c_+$ and away from $|\xi|= \lambda/c_{\pm}$, we can integrate
by parts as in the proof of Proposition \ref{prop:rio} and get an element of
$\langle z \rangle ^{-1/2+\ep}L^2(\Real^n)$.  We need only examine the
integration near $|\xi|=\lambda/c_{\pm}$ more closely.

If $c_+<c_-$, near $|\xi|=\lambda/c_-$, $y>y_M$, we have
\begin{equation}\label{eq:nearl/c_-}
(D_y+|\xi|^2-(\lambda-i0)^2)^{-1}(\widehat{Vu}(\xi,\cdot)+\hat
{f}(\xi,\cdot))(y)=e^{-iy\sqrt{\lambda^2/c_+^2-|\xi|^2}}
(a_1(\xi)+(\lambda^2/c_-^2-|\xi|^2)^{1/2}a_2(\xi))+b(y,\xi)
\end{equation}
where $a_1,a_2\in H^l(\Real^{n-1})$, $l$ an integer with $l<J-1$,
and $b\in L^2(\Real^{n-1};\lyr^{-J+1+\ep}L^2(\Real))$.
 Since again there are no
stationary points of the phase $x\cdot \xi -y\sqrt{\lambda^2/c_+^2-|\xi|^2}$
with $z/|z|$ in the support of $\chi$ and $\xi$ in the support of
$1-\Psi$, we integrate by parts.  Since, for
any $\Psi_3\in C_c^{\infty}(\Real)$,  the Fourier transform in
$x$ of $\Psi_3(|\xi|-\lambda/c_+)(\lambda^2/c_-^2-|\xi|^2)^{-1/2}$ is
in $L^p(\Real^{n-1})\cap C^{\infty}(\Real^{n-1})$ for any $p>2$, we have
the contribution of the integration of $a_i$ terms over this region is in
$\lzr^{-1/2+\ep}L^2(\Real^n)$ for any $\ep>0$.  Clearly, the
contribution of the integration of the $b$ from (\ref{eq:nearl/c_-}) is
in $\lzr^{-J+1+\ep}L^2(\Real^n)$.

Near $|\xi|=\lambda/c_+$, as in the proof of Proposition \ref{prop:rio}, we
divide the integration into two pieces.  When $|\xi|<\lambda/c_+$,
we introduce the coordinate $t=\sqrt{\lambda^2/c_+^2-|\xi|^2}$ as in
the proof of Proposition \ref{prop:rio} and can integrate by parts
to get an element of $\langle z\rangle ^{-1/2+\ep}L^2(\Real^n)$.

For the part with $|\xi|\geq\lambda/c_+$,  we use (\ref{eq:xi>l/c_+}).
Changing coordinates as in the proof of Proposition \ref{prop:rio},
we can, in the region with $\xi_1>\delta >0$, consider the
phase to be $x_1\sqrt{t^2+|\overline{\xi}|^2-\lambda^2/c_+^2}+\overline{x}\cdot
\overline{\xi}+iyt$, using $t=\sqrt{|\xi|^2-\lambda^2/c_+^2}$ and the notation
of the proof of Proposition \ref{prop:rio}.  Again we have no
stationary points and integrating by parts gives an element of
$\lzr^{-1/2+\ep}L^2(\Real^n)$.

A similar argument gives the same result for $\frac{\partial }{\partial |z|}
\chi (1-\Psi(D_x))u$.
\end{pf}

We use the notation $\Sphere^{n-1}_{\pm}=\{\om=(\omb, \omn)\in \Sphere^{n-1}:
\pm \omn >0\}$ and will use differential operators of the following type.
\begin{defin}
We say that a differential operator $P\in \Diffr^l(\Real^n)$ if $P$ is a
differential operator of the form
$$\sum_{\pm}\sum_{|\alpha|+j\leq
l}b_j(\frac{1}{|z|})a_{j,\alpha,\pm}(\frac{z}{|z|})|z|^j
(\frac{\partial}{\partial |z|}+i\lambda/c_{\pm})^jD^{\alpha}_{z/|z|}$$
where
$b_j\in C^{\infty}([0,\infty))$,
$a_{j,\alpha,\pm}\in C^{\infty}_c(\Sphere^{n-1}_{\pm})$ and
$D^{\alpha}_{z/|z|}$
is a differential operator of order $|\alpha|$ in the $z/|z|$ variables.
\end{defin}

We shall make use of the following lemma, whose proof follows by
a straightforward computation.
\begin{lemma}\label{l:cp}
If $P\in \Diffr^l(\Real^n)$ and $b\in C^{\infty}(\Sphere^{n-1})$, then
$[\frac{\partial}{\partial |z|},P] \in |z|^{-1}\Diffr^l(\Real^n),$
$[\frac{\partial}{\partial z/|z|},P] \in \Diffr^l(\Real^n)$,
$[|z|^{-m}b(z/|z|),P] \in |z|^{-m}\Diffr^{l-1}(\Real^n),$ and
$[\Delta,P] \in |z|^{-2}\Diffr^{l+1}(\Real^n).$
\end{lemma}

In order to construct the desired approximation, we shall use the
following lemma.  In practice, when we apply this lemma, the first
term will be used in solving away the error, and the subsequent terms
will be of lower order.  In particular, in applications $h$ will
vanish faster than $g$ at infinity and so will $
\chi_{\pm} (\frac{\partial}{\partial |z|}+
\frac{i\lambda}{c_{\pm}})g$.
\begin{lemma}  \label{l:ffsa}
Suppose $(\Delta -\lambda^2/c_0^2)g=h$,
$P\in \Diffr^l(\Real^n),$ and $\chi_{\pm},v\in C^{\infty}(\Sphere^{n-1})$ with
the support of $\chi_{\pm}$ contained in $\Sphere^{n-1}_{\pm}.$  Then
\begin{multline}
(\Delta -\lambda^2/c_0^2)
\left(|z|^{-j+1}v\chi_{\pm}(z/|z|)Pg\right)=
(2j-2)|z|^{-j}v\chi_{\pm}\frac{-i\lambda}{c_{\pm}}Pg
+\\
(2j-2)|z|^{-j}v\chi_{\pm}P(\frac{\partial}{\partial |z|}+
\frac{i\lambda}{c_{\pm}})g
+ |z|^{-j-1}v\chi_{\pm}(P_{l+1} g) +|z|^{-j-1}\nabla_0(v\chi_{\pm})\cdot
\nabla_0 P g
\\ + |z|^{-j+1}v\chi_{\pm}Ph+\Delta(|z|^{-j+1}v\chi_{\pm})Pg\\
= (2j-2)|z|^{-j}v\chi_{\pm}\frac{-i\lambda}{c_{\pm}}Pg+|z|^{-j-1}P'_{l+1}g
+ |z|^{-j+1}v\chi_{\pm}Ph
\end{multline}
with $P_{l+1},P_{l+1}'\in \Diffr^{l+1}(\Real^n)$.
Here $\nabla_0$ stands for the
gradient on $\Sphere^{n-1}$.
\end{lemma}
\begin{pf} The proof follows from a straight-forward computation,
using Lemma \ref{l:cp}.
\end{pf}

If $\Psi(\xi)\in C^{\infty}(\Real^{n-1})$, we use the notation
$$(\Psi(D_x)f)(z)=
(2\pi)^{1-n}\int e^{i(x-x')\cdot \xi}\Psi(\xi)f(x',y)dx'd\xi.$$
Note that $[\Psi(D_x),\Delta]=0$ and $[\Psi(D_x),c_0(y)]=0$.

\begin{lemma}\label{l:adisjsupp}
If $\chi_1,\; \chi_2\in C_c^{\infty}(\Sphere^{n-1})$ with
$\supp \chi_1 \cap \supp \chi_2=\emptyset$,
$\Psi\in {\cal S}(\Real^{n-1})$, and
$$Au =\chi_1(z/|z|)\Psi(D_x)\chi_2(z/|z|)u,$$
then $A:\langle z \rangle ^{\alpha}
L^2(\Real^n)\rightarrow \lzr^{-\infty}L^2(\Real^n)$ for any
$\alpha \in \Real$.
\end{lemma}
\begin{pf}
The crucial observation is that if $(x,y)/|(x,y)|\in \supp \chi_1$
and $(x',y)/|(x',y)|\in \supp \chi_2$, then
$|x-x'|\geq \beta |(x,y)|$,
$|x-x'|\geq \beta |(x',y)|$ for some $\beta>0$.

The Schwartz kernel of $\Psi(D_x)$ is given by $(2\pi)^{1-n}\hat{\Psi}(x'-x)$,
where $\hat{\Psi}\in {\cal S}(\Real^{n-1})$.  Thus, for $f\in L^2(\Real^n)$,
\begin{align*}
|Af|(z)& = (2\pi)^{1-n}|\chi_1(z/|z|) \int
\hat{\Psi}(x'-x)\chi_2((x',y)/|(x',y)|)f(x',y)dx'|\\
&\leq C
|\chi_1(z/|z|) \int (1+|x-x'|)^{-m-\alpha}(1+|x-x'|)^{m+\alpha}
\hat{\Psi}(x'-x)\chi_2((x',y)/|(x',y)|)f(x',y)dx'|\\
& \leq C
\lzr^{-m} (\int (1+|x-x'|)^{2m+2\alpha}
|\hat{\Psi}(x'-x)|^2dx')^{1/2} (\int \langle (x',y)\rangle^{-2\alpha}
|f(x',y)|^2dx')^{1/2}\\
& \leq C\lzr^{-m}(\int \langle (x',y)\rangle^{-2\alpha}|f(x',y)|^2dx')^{1/2}
\end{align*}
for any $m$ (where the constant depends on $m$ and $\alpha$)
and thus it follows that $Af\in \lzr^{-\infty}L^2(\Real^n)$.
\end{pf}

\begin{lemma}\label{l:psiaction}
If $\Psi \in C_c^{\infty}(\Real^{n-1})$ and $g\in \langle z \rangle
^{-\beta}L^2(\Real^n)$, then $\Psi(D_x)g\in \lzr^{-\beta}L^2(\Real^n)$.
Suppose $D^{\alpha}_zg,\;
D^{\alpha}_zPg\in \langle z \rangle ^{-\beta}L^2(\Real^n)$ for
all $P\in \Diffr^k(\Real^n)$ and for all multiindices $\alpha$.  Then
$D^{\alpha}_zP\Psi(D_x)g\in \langle z \rangle ^{-\beta}L^2(\Real^n)$ for all
$P\in \Diffr^k(\Real^n)$
 and all multiindices $\alpha$.
\end{lemma}
\begin{pf}
 To show that if $g\in \lzr^{-\beta}L^2(\Real^n)$, so is $\Psi(D_x)g$, we
take the Fourier transform of $\Psi(D_x)g$:
$${\cal F}(\Psi(D_x)g)(\eta,\tau)=\Psi(\eta){\cal F}(g)(\eta,\tau)$$
where here we are using ${\cal F}(h)(\eta,\tau)$ to denote the Fourier
transform
of $h$
in all variables.  Then, since ${\cal F}(g)\in H^{\beta}(\Real^n)$,
 ${\cal F}(\Psi(D_x)g)\in H^{\beta}(\Real^n)$,
and $\Psi(D_x)g\in \lzr^{-\beta}L^2(\Real^n)$.

We give an indication of the proof of the remainder of the lemma.
Suppose $D^{\alpha}_zg,\;
D^{\alpha}_zPg\in \langle z \rangle ^{-\beta}L^2(\Real^n)$ for
all $P\in \Diffr^k(\Real^n)$ and for all multiindices $\alpha$, and $k\geq 1$.
For $\chi_+\in C_c^{\infty}(\Sphere^{n-1}_+)$, consider
\begin{multline*}
(2\pi)^{n-1}\chi_+(z/|z|)\left(|z|\frac{\partial}{\partial |z|}+\frac{i\lambda
|z|}{c_+}\right)
\Psi(D_x)g= \chi_+(z/|z|) \int \hat{\Psi}(x'-x)
(x'\cdot \nabla_{x'}+y\frac{\partial }{\partial y} +i\lambda
|(x',y)|/c_+)f(x',y)dx'\\
+ i\lambda/c_+ \chi_+(z/|z|) \int \hat{\Psi}(x'-x) (|(x,y)|-|(x',y)|)f(x'.y)
dx'\\
-\chi_+(z/|z|)\int \sum \widehat{D_{\xi_j}\Psi}(x'-x)\frac{\partial}
{\partial x'_j}f(x',y)dx'
\end{multline*}
Let $\tilde{\chi}\in C_c^{\infty}(\Sphere^{n-1}_+)$ be $1$ on the
support of $\chi$.  Then, using the first part of the lemma,
$$\chi_+(z/|z|) \int \hat{\Psi}(x'-x)
\tilde{\chi}((x',y)/|(x',y)|)
(x'\cdot \nabla_{x'}+y\frac{\partial }{\partial y} +i\lambda
|(x',y)|/c_+)f(x',y)dx \in \lzr^{-\beta}L^2(\Real^n)$$
since
$(x'\cdot \nabla_{x'}+y\frac{\partial }{\partial y} +i\lambda
|(x',y)|/c_+)f(x',y)\in \lzr^{-\beta}L^2(\Real^n)$.  By the same reasoning,
$$\chi_+(z/|z|)\int \sum \widehat{D_{\xi_j}\Psi}(x'-x)\frac{\partial}
{\partial x'_j}f(x',y)dx'\in \lzr^{-\beta}L^2(\Real^n).$$
Moreover, by Lemma \ref{l:adisjsupp},
\begin{multline*}
\chi_+(z/|z|) \int \hat{\Psi}(x'-x)
(1-\tilde{\chi}((x',y)/|(x',y)|))
(x'\cdot \nabla_{x'}+y\frac{\partial }{\partial y} +i\lambda
|(x',y)|/c_+)f(x',y)dx
\\  \in \lzr^{-\infty}L^2(\Real^n).
\end{multline*}

To finish, note that
$\langle x-x'\rangle ^{m}\hat{\Psi}(x'-x) (|(x,y)|-|(x',y)|)$ is a smooth,
bounded, function of $x, \; y$, and $x'$ for any $m$.  Then, if $\beta \geq 0$,
\begin{align*}
& \lzr^{2\beta}|\int \hat{\Psi}(x'-x) (|(x,y)|-|(x',y)|)f(x',y)
dx'|^2  \\ &
\leq C \left| \int
\langle (x',y)\rangle^{\beta}\langle x-x'\rangle ^{\beta}
\hat{\Psi}(x'-x) (|(x,y)|-|(x',y)|)f(x',y)
dx'\right|^2 \\
& \leq C \left(\int \langle (x',y)\rangle^{2\beta}|\hat{\Psi}(x-x')
(|(x,y)|-|(x',y)|)|
|f(x',y)|^2 dx'\right) \left( \int \langle x-x'\rangle ^{2\beta}
|\hat{\Psi}(x'-x) (|(x,y)|-|(x',y)|)|dx'\right).
\end{align*}
Since $\int \langle x-x'\rangle ^{2\beta}
|\hat{\Psi}(x'-x) (|(x,y)|-|(x',y)|)|dx'<C$, where we allow the constant
$C$ to change from line to line, we have
\begin{align*}
& \int \lzr^{2\beta}|\int \hat{\Psi}(x'-x) (|(x,y)|-|(x',y)|)f(x',y)
dx'|^2 dz \\
& \leq C \int \int \langle (x',y)\rangle^{2\beta}|\hat{\Psi}(x'-x)
(|(x,y)|-|(x',y)|)|
|f(x',y)|^2 dx' dz \\
& \leq C
\int\int  \langle (x',y)\rangle^{2\beta} \left( \int(|\hat{\Psi}(x'-x)
(|(x,y)|-|(x',y)|)|) dx\right)
|f(x',y)|^2 dx' dy\\
& \leq C \int\int \langle (x',y)\rangle^{2\beta}|f(x',y)|^2dx'dy
\end{align*}
where for the last inequality we used that
$
|\int(|\hat{\Psi}(x'-x)
(|(x,y)|-|(x',y)|)|) dx| <C$.  A similar argument can be used when
$\beta<0$, using instead in the first step that for $\beta<0$,
$\langle z \rangle^{\beta} \leq C \langle x-x'\rangle^{-\beta}\langle
(x',y)\rangle^{\beta}$.

A similar argument works for a derivative in the $z/|z|$ direction, and the
argument can be iterated to get the lemma.
\end{pf}

\begin{lemma}\label{l:dr} Let $\phi, \phi_1\in C_c^{\infty}(\Real)$, with
$(1-\phi)(1-\phi_1)=1-\phi_1$ and $\phi_1(y)=1$ if $|y|\leq y_M+1$.
If $D^{\alpha}_z (1-\phi(y))g \in \langle z \rangle ^{-\beta}L^2(\Real^n)$
for $|\alpha|\leq l <\beta-1/2$, $\phi(y)g \in \lzr^{-\beta}L^2(\Real^n)$,
$\chi \in C^{\infty}_c(\Spc^{n-1})$, and
$\Psi \in C^{\infty}_b(\Real^{n-1})$ with $\supp \Psi \cap \{\xi: |\xi|=
|\lambda|/c_{\pm}\}=\emptyset$, then
$$P\chi(z/|z|)(1-\phi_1(y))
\Psi(D_x)(\Delta -(\lambda-i0)^2c_0^{-2})^{-1}g \in \langle
z \rangle ^{1/2+\epsilon}L^2(\Real^n)$$ for every $\epsilon >0$ and
$P\in \Diffr^l$, with $l<\beta-1/2$.
\end{lemma}

\begin{pf}
We give the proof for $\chi$ supported in $y/|z|>0$; the proof for $\chi$
with support in $y/|z|<0$ is similar.

We use a cut off-function, $\Psi_1\in C_c^{\infty}(\Real)$ with
$\Psi_1(|\xi|) \equiv 1 $ when $|\xi|\leq \lambda/c_+$ and supported in
a small neighborhood of that region, so that $\supp(\Psi_1(|\xi|)\Psi(\xi))
\subset \{|\xi|<\lambda/c_+\}$.
We write
\begin{multline}
(1-\phi_1(y))\chi(z/|z|)\Psi(D_x)(\Delta -(\lambda-i0)^2c_0^{-2})^{-1}g \\
= (2\pi)^{1-n}(1-\phi_1(y))\chi(z/|z|)\int e^{i x \cdot \xi}
\Psi(\xi)(\Psi_1(|\xi|)+1-\Psi_1(|\xi|))\\
(D_y^2+|\xi|^2-(\lambda-i0)^2c_0^{-2})^{-1} \hat{g}(\xi,\cdot)(y)d\xi.
\end{multline}

The main contribution is
\begin{equation}
\label{eq:mp}
(2\pi)^{1-n}(1-\phi_1(y))\chi(z/|z|)\int e^{ix\cdot \xi}
\Psi(\xi)\Psi_1(|\xi|)e^{ix\cdot \xi}
(D_y^2+|\xi|^2-(\lambda-i0)^2c_0^{-2})^{-1}\hat{g}(\xi,\cdot)(y)d\xi.
\end{equation}
Here we may write, for $y>y_M$ and $|\xi|< \lambda/c_+$,
\begin{equation}\label{eq:asy}
(D_y^2+|\xi|^2-c_0^{-2}(\lambda-i0)^2)^{-1}\hat{g}(\xi,\cdot)(y)
= e^{-iy \sqrt{\lambda^2/c_+^2 -|\xi|^2}}\tilde{g}(\xi)
+g_1(\xi,y)
\end{equation}
where $D^k_yg_1\in L^2(\Real^{n-1}_{\xi};\lyr^{-\beta+1}L^2(\Real_y))$
for all $k$.
Putting the first term of (\ref{eq:asy}) into (\ref{eq:mp}), we obtain
$$
(2\pi)^{1-n}(1-\phi_1(y))\chi(z/|z|)\int \Psi(\xi)\Psi_1(|\xi|)e^{ix\cdot \xi}
e^{-iy\sqrt{\lambda^2/c_+^2 -|\xi|^2}}
\tilde{g}(\xi)d\xi.
$$
Note that if we apply $\frac{\partial}{\partial |z|}+i\lambda/c_+$,
then the integrand vanishes on the critical set of
the phase function.  Since on the support of
$\Psi\Psi_1$, $\tilde{g}\in H^{l}$, $l$ an
integer with $l<\beta-1/2$, we
can integrate by parts to see that we have an element of
$\langle z \rangle ^{-1/2+\epsilon}L^2$.

For the tangential derivatives (in the $z/|z|$ directions),
notice that if we have a
derivative in a direction orthogonal to $y$, it commutes with
$((D_y^2+|\xi|^2)-c_0^{-2}(\lambda-i0)^2)^{-1}$.  That is, for example,
\begin{multline*}
  (-x_2\frac{\partial}{\partial x_1}+x_1\frac{\partial}{\partial x_2})
\int e^{ix\cdot \xi}\left(D_y^2+|\xi|^2-(\lambda-i0)^2c_0^{-2})^{-1}
\hat{g}(\xi,\cdot)\right)(y)d\xi\\
=
\int e^{ix\cdot \xi}\left(D_y^2+|\xi|^2-(\lambda-i0)^2c_0^{-2})^{-1}
(\xi_1\frac{\partial}
{\partial \xi_2}-\xi_2\frac{\partial}{\partial
\xi_1})\hat{g}(\xi,\cdot)\right)(y)d\xi.
\end{multline*}
By the decay properties of $g$ and the regularity and decay properties of
$(1-\phi)g$, this gives an element of $\lzr^{1/2+\ep}L^2(\Real^n)$
after multiplication by $\chi(1-\phi_1)$, if $\beta >3/2$.

 After applying
a derivative of
the form $y\frac{\partial }{\partial x_j}-x_j \frac{\partial }{\partial y}$,
as in the radial case
 the integrand vanishes on the critical set of the phase
function, and so we can integrate by parts.

This argument can be iterated up to $l<\beta -1/2$.

The second term of (\ref{eq:asy}) gives an element of
$\langle z \rangle ^{-\beta+l+1}L^2(\Real^n)$, where
$l$ is the order of the derivative
$P$.

On the support of $(1-\Psi_1(|\xi|))$,
$((D_y^2+|\xi|^2)-(\lambda-i0)^2c_0^{-2})^{-1}$ is a smooth function of
$|\xi|$, except near a finite number of points for which $\lambda^2 $ is
an eigenvalue of $c^2_0(D_y^2+|\xi|^2)$.  Since the eigenfunctions of this
operator are exponentially decreasing in $y$ and $\chi$ is supported in
$y/|z|>\delta>0$ for some $\delta >0$, these eigenfunctions do not contribute
to the asymptotics here.  Projecting off the eigenfunctions, we have
$$
(1-\Psi_1)(D_y^2+|\xi|^2-(\lambda-i0)^2c_0^{-2})^{-1}\Pi_e:
\langle y \rangle ^{-\beta}L^2(\Real_y)\rightarrow
\langle y \rangle ^{-\beta}L^2(\Real_y)
$$
with bound $C(|\xi|^2-C)^{-2}$.  Therefore,
$$
P(1-\phi_1(y))\chi(z/|z|)\int \Psi(\xi)(1-\Psi_1(|\xi|))
(D_y^2+|\xi|^2-(\lambda-i0)^2c_0^{-2})^{-1}\hat{g}(\xi,\cdot)(y)d\xi
\in \langle z \rangle ^{-\beta +l}L^2(\Real^n)$$
where we used the fact that the inverse Fourier transform is an isomorphism
on $L^2$, and the regularity properties of $D^{\alpha}_z(1-\phi(y))g$.
\end{pf}

\begin{lemma}\label{l:ps1}  Let $f\in \lzr^{-\infty}L^2(\Real^n)$,
$(1-\phi(y))f\in {\cal S}(\Real^n)$ for some $\phi \in C^{\infty}_c(\Real)$,
and let
$u=(\Delta -(\lambda -i0)^2)c^{-2})f$.
If $\Psi_0 \in C^{\infty}_b(\Real)$ is
$0$ in a neighborhood of $|\xi|=\lambda/c_{\pm}$
 and $\chi_0\in C_c^{\infty}(\Spc)$, then
there is a $w_0=\sum_{\pm} \tilde{\chi}_{\pm}\sum_{j=0}^{J-3}|z|^{-j-J+1}
P_j \Psi_0(D_x)u$, with $P_j \in \Diffr^j(\Real^{n})$, such that
$$(\Delta -\lambda^2/c_0^2)w_0=\chi_0(z/|z|) V \Psi_0(D_x)u + e_{0}$$
where $e_0 \in \lzr^{-2J+5/2+\ep}L^2(\Real^n)$.
\end{lemma}

\begin{pf}
Let $\chi_0 = \chi_+ +\chi_-$, with $\chi_{\pm}$ supported in $\pm y >0$.  We
will outline the proof for $\chi_0=\chi_+$, as the proof for
$\chi_0=\chi_-$ is similar, and the functions can be added to get the general
case.

Recall that on the support of $\chi_+$,
$V\sim \sum_{j\geq J}|z|^{-j}v_j(\frac{z}{|z|})$.  We find $w=\sum_{j=0}^{J-3}
w_{0j}$.

Let
$$w_{00}=(1-\phi_1(y))|z|^{-J+1}v_J\frac{ic_+}{\lambda(2J-2)}\chi_+\Psi_0(D_x)u$$
with $\phi_1\in C_c^{\infty}(\Real)$, $\phi_1\phi=\phi$, and $\phi_1(y)=1$
for $|y|\leq y_M$.
Then, by Lemma \ref{l:ffsa}
\begin{equation*}
(\Delta -\lambda^2/c_0^2)w_{00}=|z|^{-J}v_J(\frac{z}{|z|})\chi_+\Psi_0(D_x)u
+|z|^{-J-1}P'_1\Psi_0(D_x)u + |z|^{-J+1}\Psi_0(D_x)(Vu+f)+e_t
\end{equation*}
where $P'_1\in \Diffr^1(\Real^n)$
and $e_t\in \langle z \rangle ^{-\infty}L^2(\Real^n)$, so that
$(\Delta -\lambda^2/c_0^2)w_{00}=\chi_+V\Psi(D_x)u+e_{00}$,
$e_{00}\in \lzr^{-J-1/2+\ep}L^2(\Real^n)$.

 In the same manner, we choose $w_{01}$ so
that
\begin{equation*}
(\Delta -\lambda^2/c_0^2)w_{01}=|z|^{-J-1}v_{J+1}(\frac{z}{|z|})
\chi_+\Psi_0(D_x)u
-|z|^{-J-1}P'_1\Psi_0(D_x)u  +
|z|^{-J-2}P'_2\Psi_0(D_x)u+ e_{01r}
\end{equation*}
with $P'_2\in \Diffr^2(\Real^n)$ and $e_{01r}\in
\lzr^{-2J+1/2+\ep}L^2(\Real^n)$.
Then $(\Delta -\lambda^2/c_0^2)(w_{00}+w_{01})=V\chi_+ \Psi(D_x)+e_{01}$,
$e_{01}\in \lzr^{-J-3/2+\ep}L^2(\Real^n)$.
This can be continued, with $w_{0j}$ removing the terms in $\langle z \rangle
^{-J+1/2+\ep-j}L^2(\Real^n)$,
modulo terms in $\lzr^{-J-1/2+\ep-j}L^2(\Real^n)$,
up to $j=J-3$.
\end{pf}

We shall also need the following lemma.
\begin{lemma}\label{l:gft}
If $w(x,y)\in \langle z \rangle ^{1/2+\epsilon}L^2(\Real^n)$,
$\Psi        \in C^{\infty}_b(\Real^{n-1})$, and
$\supp\hat{w}(\xi,y) \cap \supp \Psi = \emptyset$, then
$\Psi \widehat{Vw}(\xi,y)\in H^{\infty}(\Real^{n-1}_{\xi}, \langle
y \rangle ^{-\infty}L^2(\Real_y))$.
\end{lemma}
\begin{pf}
Observe that
\begin{align}
\Psi(\xi) \widehat{Vw}(\xi,y)& = \Psi(\xi)
\int\int e^{-ix\cdot \xi}V(x,y)e^{ix\cdot \eta}\hat{w}(\eta, y)d\eta dx\\
& = \sum_j\Psi(\xi)\int\int{\Psi_j}(\xi,\eta)(\xi_j-\eta_j)^{-1}
 e^{-ix\cdot (\xi-\eta)}D_{x_j}V(x,y)
\hat{w}(\eta, y)d\eta dx
\end{align}
where $\Psi_j$ is a partition of unity with $\xi_j \not = \eta_j$ on
$\supp \Psi_j$.  We may repeat this integration by parts as many
times as desired.  Since $|D^{\alpha}_xV(x,y)|\leq C_{\alpha}
\langle z \rangle^{-J-|\alpha|}$, the lemma follows.
\end{pf}

\begin{prop}\label{prop:slot}
Suppose $\chi \in C_c^{\infty}(\Spc^{n-1})$,
 $\Psi \in C_c^{\infty}(\Real^{n-1})$ has
$\supp \Psi \cap \{ \xi:|\xi|=|\lambda|/c_{\pm}\}=\emptyset$,
$\supp \Psi \cap \{\xi: |\xi|=\kappa_j(\lambda)\}=\emptyset$,
$j=1,2,...,T(\lambda)$
, and if
$(x,y)/|(x,y)|\in supp(1-\chi)$, $y>0$, then $\pm \lambda x/(c_+|z|)\not \in
\supp \Psi$.  Moreover, suppose if $(x,y)/|(x,y)|\in supp(1-\chi)$, $y<0$, then
$\pm \lambda x/(c_-|z|)\not \in
\supp \Psi$, and
$f\in \lzr^{-\infty}L^2(\Real^{n-1})$,
$(1-\phi(y))f\in {\cal S}(\Real^n)  $ for some $\phi\in C_c^{\infty}(\Real)$.
Then
$$\chi(\frac{z}{|z|})\Psi(D_x)(\Delta -(\lambda-i0)^2/c_0^2)^{-1}f
= e^{-i\lambda|z|/c_0}\chi(\frac{z}{|z|})|z|^{-(n-1)/2}(a_0(\frac{z}{|z|})
+{\cal O}(|z|^{-1}))$$
with $\chi a_0 \in C_c^{\infty}(\Spc^{n-1})$.
\end{prop}
\begin{pf}
Recall that if $u=(\Delta -(\lambda-i0)^2/c^2)^{-1}f$, then
$u=(\Delta -(\lambda-i0)^2/c_0^2)^{-1}(Vu+f).$  If
$\Psi(\xi)\widehat{Vu+f}(\xi,y)$ were in $C^{\infty}(\Real^{n-1};
\lyr^{-\infty}L^2(\Real))$, then using
\begin{equation}\label{eq:res}
(\Delta-(\lambda-i0)^2c_0^{-2})^{-1}g=(2\pi)^{1-n}\int
e^{ix\cdot
\xi}\left((|\xi|^2+D_y^2-(\lambda-i0)^2c_0^{-2})^{-1}\hat{g}(\xi,\cdot)\right)
(y)d\xi
\end{equation}
and stationary phase, we would be done.  However, it is not clear
that $\widehat{Vu}$ should be in $C^{\infty}(\Real^{n-1};
\lyr^{-\infty}L^2(\Real))$.  We shall show that $\chi\Psi(D_x)u$ can be written
as a sum of two terms: one vanishing faster than $u$ at infinity, and
another of the form $(\Delta-(\lambda-i0)^2/c_0^2)^{-1}g_k$,
where
$\Psi(\xi)\hat{g}_k(\xi,y)\in C^{k}(\Real^{n-1};\lyr^{-k}L^2(\Real))$, where
we can make $k$ as large as desired.  Then (\ref{eq:res}) and
stationary phase
will
finish the proof.

To do this, we follow an iterative procedure.  The first step has been done in
Lemma \ref{l:ps1}.  We will iteratively construct functions $w_l$ which
have the property that
$(\Delta-\lambda^2/c_0^2)(u-w_l)$ improves with increasing $l$ in an
appropriate sense.

Let $\Psi_0, \Psi_1, \Psi_2, ...\in C_c^{\infty}(\Real^{n-1})$ be
such that, for all $i$, $\Psi_i\equiv 1$ on the support of $\Psi$,
$\Psi_{i+1}\Psi_i=\Psi_{i+1}$, and $\Psi_i$ satisfies the support requirements
placed on $\Psi$ in the statement of the Proposition.  Let $\chi_0,\;
\chi_1,\;\chi_2,...\in C_c^{\infty}(\Spc^{n-1} )$ be such that
$\chi_0\chi=\chi$ and $\chi_{i+1}\chi_i=\chi_i, $ $i=0,1,2...$.
Let $w_0$ be the
function constructed in Lemma \ref{l:ps1} for this $\Psi_0$ and $\chi_0$.
 Using the
notation of that lemma, let
$$t_0=(1-\chi_0)V\Psi_0(D_x)u+f-e_0+V(1-\Psi_0(D_x))u.$$
Note that the first three terms are in $\lzr^{-2J+5/2+\ep}L^2(\Real^n)$,
where we use the support properties of $\chi_0$ and $\Psi_0$
along with an integration by parts argument as in the proof
of Lemma \ref{l:leadorder} to obtain the
result for the first term.  Additionally, the support of the Fourier
transform in the $x$ variables of $(1-\Psi_0(D_x))u$ is disjoint
from the support of $\Psi_1$.  Then, using Lemmas \ref{l:dr} and
\ref{l:gft}, we have
$P\chi_1\Psi_1(D_x)(\Delta -(\lambda-i0)^2/c_0^2)^{-1}t_0\in \lzr^{1/2+\ep}
L^2(\Real^n)$, for all $P\in \Diffr^l(\Real^n)$, $l\leq 2J-4$.  Since
$u=w_0+(\Delta -(\lambda-i0)^2/c_0^2)^{-1}t_0$,
which can be seen by a modification of the
proof of Proposition \ref{prop:rio} and
the uniqueness result (Proposition \ref{p:powd}), this in turn means that
$P\Psi_1(D_x)w_0\in \lzr^{-J+3/2+\ep}L^2(\Real^n)$ for all $P\in \Diffr^l,$
$l\leq 2J-4$, using Lemma \ref{l:psiaction}.

Given $w_0, t_0$ as above, we now iteratively construct $w_l$ for $l\geq 1$
such that
$$(\Delta -\lambda^2/c_0^2)w_l= V\chi_l\Psi_l(D_x)u+e_l=
V\chi_l\Psi_l(D_x)w_{l-1} + V\chi_l \Psi_l(D_x)
(\Delta -(\lambda-i0)^2/c_0^2)^{-1}t_{l-1}+e_l,$$
where $e_l\in \lzr^{-J+1/2+\ep-l(J-2)}L^2(\Real^n)$ and $t_l$ is defined by
\begin{multline*}
t_l=V
\left( (1-\chi_l)\Psi_l(D_x)+
(1-\Psi_l(D_x))\right)w_{l-1}\\
+V\left((1-\Psi_l(D_x))+(1-\chi_l)\Psi_l(D_x)\right)
(\Delta -(\lambda-i0)^2/c_0^2)^{-1}t_{l-1}-e_l+f.
\end{multline*}
Then $$u=w_l+(\Delta -(\lambda-i0)^2/c_0^2)^{-1}t_l.$$
Moreover, $t_l=t_l'+t_l''$, where
$\Psi_{l+1}(\xi)\hat{t_l'}(\xi,y)\in
C^{\infty}(\Real_\xi^{n-1},\lyr^{-\infty}L^2(\Real))$,
$t_l''\in \lzr^{-(l+1)(J-2)-J+1/2+\ep}L^2(\Real^n)$,
and $P\Psi_{l+1}(D_x)w_l\in \lzr^{-J+3/2+\ep}L^2(\Real^n)$ for
all $P\in \Diffr^m(\Real^n)$, $m\leq (l+2)(J-2)$.  Additionally,
$\supp w_l(z) \subset \supp \chi_l(z/|z|)$.

Supposing that $w_{l-1},t_{l-1}$ are as above, we show how to construct
$w_l$.  Since
\begin{equation*}
P\chi_l\Psi_l(D_x)w_{l-1}\in \lzr^{-J+3/2+\ep}L^2(\Real^n) \;
\text{and} \;
 P\chi_l \Psi_l(D_x)(\Delta -(\lambda-i0)^2/c_0)^{-1}t_{l-1}\in
\lzr^{1/2+\ep}L^2(\Real^n)
\end{equation*}
 for all $P\in \Diffr^{(l+1)(J-2)}$, we can, just as
in
Lemma \ref{l:ps1}, find $w_l=\sum_{j=0}^{(l+1)(J-2)-1}w_{lj}$ so that
$$(\Delta -\lambda^2/c_0^2)w_l=V\chi_l\Psi_l(D_x)w_{l-1}+V\chi_l\Psi_l(D_x)
(\Delta -(\lambda-i0)^2/c_0^2)^{-1}t_{l-1}+e_l$$
with
$e_l\in \lzr^{-J+1/2+\ep -(l+1)(J-2)}L^2(\Real^n)$.
Let $t_l '=V(1-\Psi_l(D_x))(w_l+(\Delta-(\lambda-i0)^2)^{-1}t_{l-1}))$;
then $\Psi_{l+1}(\xi)\hat{t_l'}(\xi,y)\in
C^{\infty}(\Real_{\xi}^{n-1},\lyr^{-\infty}L^2(\Real))$ by Lemma \ref{l:gft}.
Since $\supp w_{l-1}(z)\subset \supp \chi_{l-1}(z/|z|)$,
by Lemma \ref{l:adisjsupp},
$V(1-\chi_l)\Psi_l(D_x)w_{l-1}\in \lzr^{-\infty}L^2(\Real^n)$, and
$$V(1-\chi_l)\Psi_l(D_x)(\Delta -(\lambda -i0)^2c_0^{-2})^{-1}t_{l-1}
\in \lzr^{-(l+1)(J-2)-1/2-J+\ep}L^2(\Real^n) $$ using the support properties of
$\Psi_l$ and $\chi_l$ and an integration by parts argument as in Lemma
\ref{l:dr}.  Thus $t_l''=t_l-t_l'\in \lzr^{-(l+1)(J-2)-J+1/2+\ep}
L^2(\Real^n).$

 Note that
\begin{equation}\label{eq:wlj}
w_{lj}=|z|^{-j-J+1}P_{lj}(1-\phi(y))\Psi_l(D_x)[w_{l-1}+
(\Delta-(\lambda-i0)^2/c_0^2)^{-1}t_{l-1}]=
|z|^{-j-J+1}P_{lj}(1-\phi(y))\Psi_l(D_x)u,
\end{equation}
 with $P_{lj}\in \Diffr^j(\Real^n),$ $\phi \in C_c^{\infty}(\Real)$.
Since $(\Delta -\lambda^2/c_0^2)(u-w_l)=t_l$, we have
\begin{equation}\label{eq:uis}
u=w_l+(\Delta -(\lambda-i0)^2/c_0^2)^{-1}t_l.\end{equation}

Note that
$P\Psi_{l+1}(D_x)
(\Delta -(\lambda-i0)^2/c_0^2)^{-1}t_l\in \lzr^{1/2+\ep}L^2(\Real^n)$
for $P\in \Diffr^{(l+2)(J-2)}(\Real^n)$.  Using (\ref{eq:wlj}), (\ref{eq:uis}),
and Lemma \ref{l:psiaction}, this in turn means that $P\Psi_{l+1}(D_x)w_l
\in \lzr^{3/2+\ep -J}L^2(\Real^n)$ for all $P\in \Diffr^{(l+2)(J-2)}(\Real^n)$.
Thus, for any $l\geq 1$, $w_l$ and $t_l$ can be constructed to have the
desired properties.

To prove the proposition, we use (\ref{eq:uis}).  Since
$t_l=t_l'+t_l''$ with $\Psi_{l+1}(\xi)\hat{t_l'}(\xi,y)\in
C^{\infty}(\Real_\xi^{n-1},\lyr^{-\infty}L^2(\Real))$,
$t_l''\in \lzr^{-(l+1)(J-2)-J+1/2+\ep}L^2(\Real^n)$, we have
$\Psi(\xi)\hat{t_l}(\xi,y)\in H^s(\Real^{n-1};\lyr^{-(l+1)(J-2)+1/2+\ep-J+s}L^2
(\Real))$
for $s<(l+1)(J-2)+J-1/2-\ep$.  Then, using equation (\ref{eq:res}) and
stationary
phase, we see that
$$\chi\Psi(D_x)(\Delta -(\lambda-i0)^2/c_0^2)^{-1}t_l
=
\chi e^{-i\lambda |z|/c}|z|^{-(n-1)/2}(a_0(\frac{z}{|z|})
+{\cal O}(|z|^{-1}))$$
with $a_0\in C^{(l+1)(J-2)+J-n-4}(\Spc)$ when $l$ is sufficiently large.

To finish, then, we only need show that $\Psi(D_x)w_l$ is of order
$|z|^{-(n+1)/2}$.
But we recall that $$D^{\alpha}_{z/|z|}w_l,\;
|z| (\frac{\partial}{\partial |z|}+i\lambda/c)D^{\alpha}_{z/|z|}w_l
\in \lzr^{-J+3/2+\ep}L^2(\Real^n).$$  This in turn means that
$\Psi(D_x)w_l$ has the same properties, and that $\Psi(D_x)w_l
={\cal O}(|z|^{-(n-1)/2-J+1+\ep})$.
\end{pf}

\subsection{The structure of the scattering matrix}\label{s:sosm}

Combining the construction of the
approximate Poisson operator
of Section \ref{s:cpo} and Theorem \ref{thm:smc}, we
have now proved Theorem \ref{thm:sosm}.  For completeness, we remark that
the construction of the approximation
of $P(\lambda)\Pi_j$, $1\leq j \leq T(\lambda)$
in Section \ref{s:pj}, Corollary \ref{c:evp}, and Theorem \ref{thm:smc} prove
\begin{prop}\label{p:notmain}
Let $c,c_0$ satisfy the general conditions of Section \ref{s:aan} and either
hypothesis (H1) or (H2).  Let $A(\lambda)=(A_{ij}(\lambda))$,
$0\leq i,j\leq T(\lambda)$.  Then, for $j>0$, $A_{jj}(\lambda)$ is
a Fourier integral operator associated with the antipodal mapping on
$\Sphere^{n-2}$, and $A_{ij}(\lambda)$ is a smoothing operator when
$i\not = j$.
\end{prop}

\section{The Inverse Problem}\label{s:inverse}

We recall our central inverse result, Theorem \ref{thm:inverse}:
\begin{thmb}
Suppose $c$ and $c_0$ satisfy the general
assumptions of Section \ref{s:aan}, as well as either hypothesis (H1) or (H2),
and $n\geq 3$.  Then,
if $c_+=c_-$, the asymptotic expansion at infinity of $c-c_0$ is uniquely
 determined by $c_0$ and the transmitted singularities of the main part of the
scattering matrix at fixed nonzero
energy.  If $c_+<c_-$, then the asymptotic expansion is uniquely determined by
$c_0$ and the reflected singularities of the main part of the scattering matrix
at fixed nonzero energy.
\end{thmb}

In proving the results for the inverse problem, we use the techniques of
\cite{smagailess, recpoten} and much of their language.  We recall the
arguments from these papers, noting the adjustments that must be made
for the stratified case.

Theorem \ref{thm:inverse} follows
 from the following theorem, which is somewhat
stronger.
\begin{thm}\label{thm:sinverse}
Suppose $n\geq 3$,
$c_1$ and $c_2$ satisfy the general requirements of Section
\ref{s:aan}, and either (H1) or (H2), for the same $c_0$.
Let $S_1(\lambda)$, $S_2(\lambda)$ be the corresponding scattering
matrices for some $\lambda \in \Real \setminus \{0\}$.  If, for $c_+=c_-$,
the transmitted part of the main part of $S_1(\lambda)
-S_2(\lambda)$ is of order $-l$, then
$c_1(z)-c_2(z)= {\cal O}(|z|^{-l-1}).$
If for $c_+<c_-$ the reflected part of $S_1(\lambda)
-S_2(\lambda)$ is of order $-l$, then $c_1(z)-c_2(z)={\cal O}(|z|^{-l-1})$.
\end{thm}
\begin{pf}
Let $S^{k}_{\cl,s}(\Real^n)$ be the set of functions
which, for $|z|>1$, are of the form
$\sum_{j\geq 0}|z|^{k-j}a_{j}(z/|z|)$, where $a_{j}(z/|z|)\in C^{\infty}_b
(\Sphere^{n-1}\setminus\{(\omb,0)\}).$

Suppose that $c_1-c_{2}=W\in S^{-k}_{\cl,s}(\Real^n)$ and that the
scattering matrices associated to
$c_1^2\Delta$ and $c_2^2\Delta$ have the same transmitted (if $c_+=c_-$)
or reflected (if $c_+\not = c_-$) singular parts, to
order $l\geq k$.  Then we shall show that
actually $c_1-c_{2}\in S^{-k-1}_{\cl,s}(\Real^n)$, and thus by
induction $c_1-c_{2}\in S^{-l-1}_{\cl,s}(\Real^n).$

If $c_1-c_2\in S^{-k}_{\cl,s}(\Real^n)$, then
$\lambda^2(c^{-2}_1-c^{-2}_2)=\lambda^2
c_1^{-2}c_2^{-2}(c_1+c_2)(c_2-c_1)=U$,
with
\begin{equation*}
U_{\rest y>y_M}\sim \sum_{j\geq k}|z|^{-j}U_{-j,+}(z/|z|),\;\;
U_{\rest y<-y_M}\sim \sum_{j\geq k}|z|^{-j}U_{-j,-}(z/|z|),
\end{equation*}
and $U={\cal O}(|z|^{-k})$, where $U_{-j, \pm}\in C^{\infty}_b
(\Sphere^{n-1}_{\pm})$.    Let \begin{equation*}
W_{-j}(\omb,\omn)=\left\{
\begin{array}{cc} U_{-j,+}(\omb, \omn) & \text{if}\; \omn>0\\
U_{-j,-}(\omb, \omn) & \text{if}\; \omn<0
\end{array}\right.
\end{equation*}
Note then that
the first
$k-2$ terms in the construction of the Poisson operator carried out in
Section \ref{s:cpo} are the same, and the difference comes in the $k-1$st
term.

Although many of the underlying techniques are the same, we shall treat the
cases $c_+=c_-$ and $c_+\not = c_-$ serially.

We begin with the case $c_+=c_-$.
In the construction of the Poisson operator, the transmitted parts (that
is, in the lower hemisphere)
of the $k-1$st terms differ by
\begin{equation}
\frac{i |z|^{1-k}}{2\lambda c_+ (\sin s)^{k-1}}
T_+(\omn, \lambda) \int _0^s W_{-k}(s',\theta;\om)(\sin s')^{k-2}ds',
\end{equation}
almost as in (2.3) of \cite{smagailess}.  Here $T_+(\lambda, \omn)$ is the
transmission coefficient determined by equations (\ref{eq:ode})-(\ref{eq:yni}).
We remark that in case $c_+=c_-$, $T_+$=$T_-$.  The transmission coefficient
must be nonzero for $\lambda \in \Real$, $\lambda \not =0$, $\omn \not =0$,
$0<\omn<1$.  Therefore, as described in \cite[Section 4]{smagailess}, we can
recover from the difference of the transmitted parts of the scattering
matrices
$$\int _0^{\pi}W_{-k}(s,\theta;\om)(\sin s)^{k-2}ds$$
as long as $\om=(\omb,\omn)$ with $\omn \not = 0$.

If the transmitted parts of the two scattering matrices are the same to
order $k-1$, then \begin{equation}
I_k=\int_0^{\pi}W_{-k}(s,\th;\om)(\sin s)^{k-2}ds=0.
\end{equation}
Since this is true for all $\om$ with $\omn \not = 0$, we can differentiate
with respect to the starting point twice, use
$\sin^2 s+ \cos^2 s=1$, and find that $I_{k-2}=\int_0^{\pi}W_{-k}(s,\theta;\om)
(\sin s)^{k-4}ds=0$
as well (see \cite[Section 5]{recpoten}).
Therefore, if $k$ is even, we reduce eventually to the case with $k-2=0$ and
if $k$ is odd, to $k-2=1$.  When
$k$ is even, differentiating one more time with respect to the starting
point shows that $W_{-k}$ is even; for odd $k$ two more differentiations
show that $W_{-k}$ is odd.

When $k$ is even, we obtain
$$\int_{\gamma}W_{-k}=0$$
for each closed geodesic $\gamma$ by joining together two half-geodesics.
However, by \cite[Theorem 4.7]{helgason}, for $n\geq 3$ the $x$-ray transform
on $\Sphere^{n-1}$ with domain restricted to smooth
even functions is 1-1.  Although $W_{-k}$ may have a jump discontinuity at
$\omn=0$, it is smooth elsewhere.  As in the proof of
\cite[Corollary 4.19]{helgason}, by first taking a convolution with
$W_{-k}$ we may assume that $W_{-k}$ is smooth and, applying the
theorem, it is thus $0$.

If $k$ is odd, we consider $\frac{z_n}{|z|}W_{-k}$, which is even.
Since for each geodesic beginning at $z_n=0$, $\sin s$ is a constant multiple
of $z_n/|z|$, we obtain $\int_{\gamma}\frac{z_n}{|z|}W_{-k}=0$ for
each geodesic $\gamma$, and again $W_{-k}=0$.

When $c_+\not = c_-$, we use the reflected singularities in the inverse
problem.  Recalling that $W\in S^{-k}_{\cl,s}$,
we can recover from the reflected singularities, when $\omn>0$,
\begin{equation}\label{eq:rpart}
R_+(\omn, \lambda)\left( \int _0^{s_0}W_{-k}(s,\th;\om)(\sin s)^{k-2}ds
+ \int _{s_0}^{\pi}W_{-k}(s',\th;\om)(\sin s')^{k-2}ds'.
\right)
\end{equation}
The first integral is along a geodesic originating at $\om$ and continuing
to $\{(\overline{\phi},0)\}\subset \Sphere^{n-1}$; the second integral
is along the reflection of the first geodesic when it meets $\phi_n=0$
and the path of integration ends at the point $(-\omb, \omn)$.  The variable
$s'$ is the distance from the point $(\omb,-\omn)$.

It is, however, more convenient to think of the sum (\ref{eq:rpart})
as the single integral
\begin{equation}\label{eq:nicer}
R_+(\omn, \lambda) \int _0^{\pi}\tilde{W}_{-k+}(s,\th;\om)(\sin s)^{k-2}ds
\end{equation}
where
\begin{equation}
\tilde{W}_{-k+}(\phi)=\left\{ \begin{array}{ll}
 W_{-k}(\phi) & \text{if} \;\phi_n\geq 0\\
W_{-k}(\overline{\phi}, -\phi_n) & \text{if} \;\phi_n<0
\end{array}
\right.
\end{equation}
and $s$ is the distance from $\om$.  It is fairly straightforward to see by
symmetry that (\ref{eq:rpart}) and (\ref{eq:nicer}) are the same.

If we can show that (\ref{eq:nicer}) is sufficient for recovering
$\int_0^{\pi}\tilde{W}_{-k,+}(s,\theta;\om)(\sin s)^{k-2}ds$ for
all $\om$ with $\omn>0$, then the analysis used in the case $c_+=c_-$
will show that $\tilde{W}_{-k,+}$ is $0$ if the reflection coefficients
agree to order $-k$.

It suffices that
 $R_+(\omn,\lambda)$ is $0$ for at most an isolated set of $\omn$ with
$0<\omn \leq 1$, for we can obtain the integrals for these isolated values
of $\omn$ by continuity.  We recall that for $0<\omn<\sqrt{1-c_+^2/c_-^2}$,
$|R_{+}(\omn,\lambda)|=1$ (\cite[Chapter 3]{wilcox}).
Moreover, because $c_0(y)-c_{\pm}
$ is compactly supported for $\pm y >0$,
for fixed $\lambda \in \Real$, $R_{+}(\omn,\lambda)$ can be
extended to a meromorphic function of $\omn$ in a neighborhood of $0<\omn <1$,
except near $\omn =\sqrt{1-c_+^2/c_-^2}$, where it is a meromorphic function
of $(1-c_-^2/c_+^2+c_-^2\omn^2/c_+^2)^{1/2}$.  Therefore, its zeros are
 isolated, and we have shown that it is possible to recover $W_{-k}(\om)$
for $\omn>0$.

A very similar analysis works for the lower hemisphere, proving the theorem.
\end{pf}

We remark that this proof shows that if $c_1-c_2\in S_{\cl}(\Sphere^{n-1})$,
then the main part of
$S_1(\lambda)-S_2(\lambda)$ is of order $-k+1$.

\begin{cor}\label{c:inverse}
Let $c$, $c_0$ satisfy the general conditions of Section \ref{s:aan}, and
either (H1) or (H2), and let $n \geq 3$.
Then $c_+$, $c_-$, and the main part of the scattering
matrix at nonzero fixed energy determine $c$ modulo terms vanishing faster
than the reciprocal of any polynomial at infinity.
\end{cor}
\begin{pf}
We need only show that $c_0$ can be recovered from $c_+$, $c_-$, and
knowledge of the scattering matrix at fixed energy.  The leading order
singularity of the scattering matrix $A(\lambda)$ determines and is determined
by $R_{\pm}(\lambda, \omn)$, $T_{\pm}(\lambda, \omn)$, $\lambda$
and $c_{\pm}$, where $R_{+}$, $T_+$
are defined by equations (\ref{eq:ode}-\ref{eq:yni}), and a similar definition
gives $R_-$, $T_-$, just as in one-dimensional scattering theory (see e.g.
\cite{c-k}).

Fix $\lambda \in \Real \setminus \{0\}$.
We can think of equation (\ref{eq:ode}) in the slightly more
general form
\begin{equation}\label{eq:seq}
 (D_y^2-\lambda^2(1/c_0^2-1/c_+^2)-k^2)\phi=0,
\end{equation}
a Schr\"odinger operator with potential $-\lambda^2(1/c_0^2-1/c_+^2)$ which
is either compactly supported (if $c_+=c_-$) or ``steplike'' (if
$c_+<c_-$).  We can define the reflection and transmission coefficients
$\tilde{R}_{\pm}(k)$, $\tilde{T}_{\pm}(k)$,
for (\ref{eq:seq}) as usual, as in (\ref{eq:ode}-\ref{eq:yni}),
and $\tilde{R}_{\pm}(k)=R_{\pm}(\lambda,\omn)$, $\tilde{T}_{\pm}(k)
=T_{\pm}(\lambda, \omn)$, when $k=\lambda \omn/c_+$.  Moreover,
$\tilde{R}_{\pm}$, $\tilde{T}_{\pm}$ are meromorphic functions of
$k\in \Complex$ if $c_+=c_-$, and if $c_+<c_-$, they are meromorphic functions
on $ \hat{Z}$, the Riemann surface on which $k$ and
$(k^2-\lambda^2/c_+^2+\lambda^2/c_-^2)^{1/2}$
are single-valued holomorphic functions.  Therefore, knowing
$R_+(\lambda, \omn)$, $T_{+}(\lambda, \omn)$ for
$0<\omn<\sqrt{1/c_+^2-1/c_-^2}$ determines $\tilde{R}_+(k)$, $\tilde{T}_+(k)$
on
the whole plane (if $c_+=c_-$) or $\hat{Z}$ (if $ c_+<c_-$).  This in
turn determines the eigenvalues of $D_y^2-\lambda^2(1/c_0^2-1/c_+^2)$ and
the norming constants.  These, together with $c_{\pm}$ and
$\tilde{R}_+$, determine the potential $-\lambda^2(1/c_0^2-1/c_+^2)$
(e.g. \cite{d-t,c-k}).
\end{pf}

\small
\noindent
{\sc
Department of Mathematics,
University of Missouri
Columbia, Missouri 65211 U.S.A.\\
e-mail: tjc@@math.missouri.edu
}

\vspace{2mm}

\noindent
{\sc Royal Bank of Scotland Group Risk,
Level 9,
135 Bishopsgate,
London,
EC2M 3UR, U.K.\\
e:mail: markjoshi@@alum.mit.edu}

\end{document}